\documentclass[a4paper]{article}
\usepackage[T1]{fontenc}
\usepackage[latin9]{inputenc}
\usepackage[verbose,hmargin=1in]{geometry}
\usepackage{mathtools}
\usepackage[svgnames]{xcolor}
\usepackage{amscd}
\usepackage{tikz-cd}
\usepackage{amsmath}
\usepackage{amsthm}
\usepackage{amssymb}
\usepackage{setspace}
\usepackage{cite}
\usepackage[bookmarks=true]{hyperref}
\usepackage{bookmark}

\numberwithin{equation}{section}
\numberwithin{figure}{section}

\theoremstyle{plain}
\newtheorem{thm}{Theorem}
\newtheorem{lem}[thm]{Lemma}
\newtheorem{prop}[thm]{Proposition}
\newtheorem{cor}[thm]{Corollary}
\newtheorem*{thmiv}{Theorem 5}
\newtheorem*{thmv}{Theorem 6}

\theoremstyle{definition}
\newtheorem{defn}[thm]{Definition}
\theoremstyle{remark}
\newtheorem{rem}{Remark}

\newcommand{\ax}{\mathrm{axial}}
\newcommand{\Dreg}{\mathcal{D}_\mathrm{reg}}

\begin{document}
\title{Maslov $S^{1}$ Bundles and Maslov Data}
\author{%
Konstantinos Efstathiou$^a$,
Bohuan Lin$^b$, 
Holger Waalkens$^c$
}
\maketitle
\begin{center}
{$^a$} Zu Chongzhi Center, Duke Kunshan University, Kunshan, Jiangsu, China. \\
Email: k.efstathiou@dukekunshan.edu.cn \\
{$^b$} School of Mathematics and Physics, Xi'an Jiaotong Liverpool University, Suzhou, Jiangsu, China.  \\
Email: bohuan.lin@xjtlu.edu.cn \\
{$^c$} Bernoulli Institute, University of Groningen, Groningen, The Netherlands. \\ 
Email: h.waalkens@rug.nl \\
\end{center}
\medskip

\begin{abstract}
We define Maslov $S^1$ bundles over a symplectic manifold $(M,\omega)$.
These are the determinant bundle $\Gamma_J$ of the unitary frame bundle defined by an almost complex structure compatible with $\omega$, and the bundle $\Gamma_J^2 = \Gamma_J \big/ \{\pm1\}$. 
We analyze the properties of the Maslov $S^1$ bundles $\Gamma_J$ and $\Gamma_J^2$, focusing on the interplay between their geometry and the dynamics of a symplectic action of a compact Lie group $G$ on $M$ which induces lifted $G$ actions on $\Gamma_J$ and on $\Gamma_J^2$.
We show that when $M$ is a homogeneous $G$-space and the first real Chern class $c_\Gamma$ is nonvanishing, $\Gamma_J$ and $\Gamma_J^2$ are also homogeneous $G$-spaces. 
Moreover, we give an alternative proof of the fact that when $[\omega]=r\,c_{\Gamma}$ for some real number $r$, then the symplectic $G$ action on $(M,\omega)$ is Hamiltonian.
When the Maslov $S^1$ bundle $\Gamma_J^2$ is trivial, then an index generalizing the Maslov index can be defined.
This is no longer true if $\Gamma_J^2$ is not trivial.
However, if $G=S^1$ acts symplectically on $(M,\omega)$ we define a quantity that we call Maslov data which serves as a non-integrable version of the notion of Maslov index in the case where $\Gamma_J^2$ is not trivial, and we associate the Maslov data at fixed points of the $G=S^1$ action to their resonance type.
Finally, we consider three applications motivated by the study of integrable Hamiltonian systems. 
First, we discuss conditions under which an $S^1$ symmetry of a two degrees of freedom integrable Hamiltonian system can be extended to a $\mathbb T^2$ symmetry. 
Second, we show that the Maslov $S^1$ bundles over Lagrangian pinched tori are trivial.
Third, we consider $S^2 \times S^2$ as a symplectic manifold with an $S^1$ action corresponding to simultaneous rotations of the two spheres, and we compute the corresponding Maslov data.
\end{abstract}

\section{Introduction}
\label{sec:intro}

Maslov indices originated in the study of semi-classical approximations in quantum mechanics \cite{Maslov1972}. 
They are playing an important role for quantization as well as giving insights into geometric aspects of Hamiltonian dynamics and symplectic manifolds \cite{Arnold1967a, Audin2014, Dullin2005, McDuff2017}.
Arnold proved in \cite{Arnold1967a} that the Maslov index of a closed path $\tilde{\gamma}$ of Lagrangian planes in $\mathbb{R}^{2n}$ can be characterized either as the (signed) number of intersections of $\tilde{\gamma}$ with the Maslov cycle $\mathfrak{M}_0$ induced by the Lagrangian subspace $\mathbb{E}_{0} = \mathbb{R}^{n}\times\{0\}$,
or, equivalently, as the degree of the Maslov-Arnold map for $\tilde{\gamma}$:
\begin{equation}
\mathfrak{m}_{\tilde{\gamma}}:S^{1}\xrightarrow{\tilde{\gamma}}\Lambda(n)\cong\mathbb{U}(n)\big/\mathbb{O}(n)\xrightarrow{\det_{\mathbb{C}}^{2}}S^{1},\label{eq:def_Maslov-Arnold}
\end{equation}
where $\Lambda(n)$ is the Lagrangian Grassmannian of $\mathbb{R}^{2n}$.

Maslov indices are defined in \cite{Maslov1972} and further analyzed in \cite{Arnold1967a} in the phase space $(\mathbb{R}^{2n}, \sum_j dx_j\wedge dy_j$).
In the case of a general symplectic manifold $(M,\omega)$, the Lagrangian Grassmannian is replaced by the bundle of Lagrangian planes $\Lambda_{pl}$, which consists of all the Lagrangian subspaces of the tangent spaces of $M$, namely, 
\begin{equation}
\Lambda_{pl} = \bigsqcup_{p\in M}\{\text{Lagrangian\ subspaces\ of\ $T_{p}M$} \},
\label{eq:Lambda_pl_over_M}
\end{equation}
see \cite{Butler2009, Contreras2003}.
In this context, a Maslov cycle $\mathfrak{M}$ consists of all the Lagrangian spaces that have nontrivial intersection with a Lagrangian subbundle $\mathfrak{L}$, and the Maslov index of a loop in $\Lambda_{pl}$ is then defined, in the same manner as in \cite{Arnold1967a}, as the number of intersections with $\mathfrak{M}$.

At the same time, the bundle of Lagrangian planes $\Lambda_{pl}$ can also be viewed as a structure associated with an almost complex structure $J$ compatible with $\omega$.
In the case $(M,\omega) = (\mathbb{R}^{2n}, \sum_j dx_j \wedge dy_j$), the bundle of Lagrangian planes is the trivial bundle
\[
\Lambda_{pl} =
\mathbb{R}^{2n}\times\Lambda(n)=\mathbb{R}^{2n}\times\mathbb{U}(n)\big/\mathbb{O}(n),
\]
where the last space can be interpreted as the quotient space of $\mathbb{R}^{2n}\times\mathbb{U}(n)$ by the $\mathbb{O}(n)$ action on $\mathbb{U}(n)$. 
Here, $\mathbb{R}^{2n}\times\mathbb{U}(n)$
is the unitary frame bundle consisting of all the unitary frames of
the tangent spaces of $\mathbb{R}^{2n}$ with respect to the canonical
almost complex structure
\[
J_0 = \begin{bmatrix}
0 & -I_n \\ I_n & 0
\end{bmatrix}
\]
on the tangent bundle $T\mathbb{R}^{2n}$. 
The $S^{1}$ appearing at the end of the sequence \eqref{eq:def_Maslov-Arnold} as the codomain of $\mathfrak{m}_{\tilde{\gamma}}$ can be thought of as the fibers of the $S^{1}$ principal bundle 
\[
\mathbb{R}^{2n}\times\mathbb{U}(n)\big/\mathbb{EU}(n)\cong\mathbb{R}^{2n}\times S^{1},
\]
where $\mathbb{EU}(n)$ is the subgroup of $\mathbb{U}(n)$ consisting
of all the matrices with complex determinant $\pm1$.

Our aim in this paper is to generalize these notions from the phase space $(\mathbb{R}^{2n}, \sum_j dx_j\wedge dy_j)$ to arbitrary symplectic manifolds $(M,\omega)$, and to analyze their properties.
Recall that for an arbitrary symplectic manifold $(M,\omega)$, an almost complex structure $J$ compatible with the symplectic structure $\omega$ can always be chosen \cite{McDuff2017, CannasDaSilva2001}. 
We can then define the counterparts of $\mathbb{R}^{2n}\times\mathbb{U}(n)$, $\mathbb{R}^{2n}\times\Lambda(n)$, and $\mathbb{R}^{2n}\times\mathbb{U}(n)\big/\mathbb{EU}(n)$ as follows.

\begin{defn}
Given a symplectic manifold $(M,\omega)$ with compatible almost complex structure $J$, the $\mathbb{U}(n)-$principal bundle 
\begin{equation}
\label{eq:unitary-frame-bundle}
\pi_{Fr_{J}^{u}}:Fr_{J}^{u}\rightarrow M
\end{equation}
associated to $J$, is called the \emph{unitary frame bundle}. 
Here, $\mathbb{U}(n)$ acts from the right and $\pi_{Fr_{J}^{u}}^{-1}(p) \cong \mathbb{U}(n)$ is the space of unitary frames of $T_pM$.
Then, the fiber bundle 
\begin{equation}\label{eq:lagrangian-plane-bundle}
\pi_{\Lambda_{J}} : 
\Lambda_{J}=Fr_{J}^{u}\big/\mathbb{O}(n) \to M,
\end{equation}
is called the \emph{Lagrangian plane bundle}.
Finally, the principal $S^1$ bundles
\begin{equation}\label{eq:maslov-s1-bundles}
\pi_{\Gamma_J^2} : 
\Gamma_{J}^{2}=Fr_{J}^{u}\big/\mathbb{EU}(n) \to M,
\quad 
\pi_{\Gamma_J} : 
\Gamma_J=Fr_{J}^{u}\big/\mathbb{SU}(n) \to M,
\end{equation}
are both called \emph{Maslov $S^1$ bundles}.
\end{defn}

\begin{rem}
The Maslov $S^1$ bundle $\Gamma_J^2$ is the counterpart of $\mathbb{R}^{2n}\times\mathbb{U}(n)\big/\mathbb{EU}(n)$, and it will be used to define the Maslov index. 
However, it will also be useful to consider the Maslov $S^1$ bundle $\Gamma_J$.
The two Maslov $S^1$ bundles are related by $\Gamma_{J}^{2} \cong \Gamma_{J}\big/\{\pm1\}$, where the quotient map $q_{\pm1} : \Gamma_J \to \Gamma_J^2$ is induced by the quotient map $\mathbb{EU}(n) \to \mathbb{SU}(n) = \mathbb{EU}(n) \big/ \{\pm1\}$.
\end{rem}

\begin{rem}
Although from the set-theoretic point of view, the definitions of $Fr_{J}^{u}$, $\Lambda_{J}$, $\Gamma_J^2$, and $\Gamma_{J}$ depend on the choice of the compatible almost complex structure $J$, their bundle structures are independent of this choice. 
\end{rem}

The bundle $\Lambda_{J}=Fr_{J}^{u}\big/\mathbb{O}(n)$ is in one-one correspondence to $\Lambda_{pl}$. 
This is because each fiber $\pi_{\Lambda_{J}}^{-1}(p)$ of $\Lambda_{J}$ with $p\in M$ is exactly the Lagrangian Grassmannian of the symplectic vector space $(T_{p}M,\omega_{p})$
with respect to the compatible linear complex structure $J\big|_{T_{p}M}$.
Moreover, as we discuss in detail in Appendix~\ref{sec:topology-lambda}, the spaces $\Lambda_J$ and $\Lambda_{pl}$ are homeomorphic with respect to their natural topologies and, for this reason, we use them interchangeably. 
From the set-theoretic point of view, the advantage of using $\Lambda_{pl}$ is that for any symplectic flow
$\varphi^t$ on $M$, the tangent maps $\varphi^t_*$ induce a
flow on $\Lambda_{pl}$ in a natural way. 
The advantage of considering $\Lambda_{J}$ is its direct relation with $\Gamma_J^2$ and $\Gamma_J$.

Since $\mathbb{O}(n)\subset\mathbb{EU}(n)$, and they are both closed subgroups of $\mathbb{U}(n)$, each $\mathbb{O}(n)$ orbit in $Fr_{J}^{u}$ lies in a single $\mathbb{EU}(n)$ orbit, and there is a natural quotient map
\begin{equation} \label{eq:det^2_J}
det_{J}^{2} : \Lambda_J = Fr_{J}^{u}\big/\mathbb{O}(n) \rightarrow \Gamma_J^2 = Fr_{J}^{u}\big/\mathbb{EU}(n),
\end{equation}
sending each $\mathbb{O}(n)$ orbit in $Fr_{J}^{u}$ to its containing $\mathbb{EU}(n)$ orbit. 
Viewed in proper local trivializations of the bundles $\Lambda_{J}=Fr_{J}^{u}\big/\mathbb{O}(n)\rightarrow M$
and $\Gamma_{J}^{2}=Fr_{J}^{u}\big/\mathbb{EU}(n)\rightarrow M$,
the quotient map $det_J^2$ takes the form
\[
U\times\mathbb{U}(n)\big/\mathbb{O}(n)\ni(p,[A])\mapsto(p,det_{\mathbb{C}}^{2}A)\in U\times S^{1}
\]
with respect to local charts over an open set $U \ni p$ of $M$.
The composition
\begin{equation} \label{eq:maslov-arnold-map}
det_{J}^{2}\circ\tilde{\gamma}:S^{1}\rightarrow\Gamma_{J}^{2}
\end{equation}
is then the Maslov-Arnold map for a closed path $\tilde{\gamma}$
in $\Lambda_{J}$. 

When the Maslov $S^1$ bundle admits a trivialization $\Gamma_{J}^{2}\cong M\times S^{1}$ with respect to some global section $\mathfrak{s}^2$, such a trivialization leads to the definition of an index.

\begin{defn}\label{def:maslov-index(intro)}
Given a global section $\mathfrak{s}^2$ of $\Gamma_J^2$, inducing a trivialization $tr_{\mathfrak{s}^2}: \Gamma_{J}^{2} \to M\times S^{1}$, the \emph{index} $\widetilde{\mathfrak{m}}_{\mathfrak{s}^2}(\tilde{\gamma})$ for a closed path $\tilde{\gamma}$ of Lagrangian planes in $\Lambda_J$ is the degree of the map
\begin{equation}
S^{1}
\xrightarrow{\tilde{\gamma}}\Lambda_{J}
\xrightarrow{det_{J}^{2}} \Gamma_{J}^{2}
\xrightarrow{tr_{\mathfrak{s}^2}} M\times S^{1}
\xrightarrow{pr_{S^{1}}} S^{1}.
\label{eq:ind_m_s_g}
\end{equation}
\end{defn}

Note that a globally defined Lagrangian vector subbundle gives a global section $\sigma: M \to \Lambda_{J}$, and then the composition $\mathfrak{s}^2 = det_{J}^{2}\circ\sigma$ is a global section to $\Gamma_J^2$, which gives a trivialization $\Gamma_J^2\cong M\times S^1$.
Also note that this is always the case if $M$ is the cotangent bundle of some manifold, and the Lagrangian subbundle can be chosen as the vertical distribution. 
In this case, the index $\widetilde{\mathfrak{m}}_{\mathfrak{s}^2}(\tilde{\gamma})$ defined above is exactly the Maslov index (with respect to $\mathfrak{s}^2$) in the usual sense.

In particular, suppose that $\mathcal{S}$ is a Lagrangian submanifold of $M$.
In this case, a loop $\gamma$ in $\mathcal{S}$ lifts to the loop $\sigma_{\mathcal{S}}\circ\gamma : t \mapsto T_{\gamma(t)}\mathcal{S}$ in $\Lambda_{J}$, 
where $\sigma_{\mathcal{S}} : \mathcal{S} \to \Lambda_J \big|_{\mathcal{S}}$ is the smooth section which assigns to each $b \in \mathcal{S}$ the tangent space $T_b\mathcal{S}$. 

\begin{defn}\label{def:maslov-index-lagrangian}
If $\mathcal S$ is a Lagrangian submanifold of $M$, and $\gamma : S^1 \to \mathcal S$ is a loop in $\mathcal S$, then the \emph{Maslov index of $\gamma$} with respect to a global section $\mathfrak{s}^2$ of $\Gamma_J^2$ is
\[ \mathfrak{m}_{\mathfrak{s}^2}(\gamma) 
= \widetilde{\mathfrak{m}}_{\mathfrak{s}^2}(\sigma_{\mathcal{S}}\circ\gamma), \]
that is, $\mathfrak{m}_{\mathfrak{s}^2}(\gamma)$ is the degree of the map
\[
S^{1}
\stackrel{\gamma}{\longrightarrow}
\mathcal{S}
\stackrel{\sigma_{\mathcal{S}}}{\longrightarrow}
\Lambda_{J}
\stackrel{det_{J}^{2}}{\longrightarrow}
\Gamma_{J}^{2}\stackrel{tr_{\mathfrak{s}^2}}{\longrightarrow}
M\times S^{1}
\stackrel{pr_{S^{1}}}{\longrightarrow}
S^{1}.
\]
\end{defn}

Generally speaking, the Maslov $S^{1}$ bundle $\Gamma_{J}^{2}$ is not necessarily trivial as, for example, when $M=S^2$ with the standard symplectic form, see Sec.~\ref{sec:example-S2}. 
In this case, the mapping \eqref{eq:ind_m_s_g} is not defined.
Recall, however, that in the standard context of the phase space $(\mathbb{R}^{2n}, \sum_j dx_{j}\wedge dy_{j})$ with (piecewise) smooth loop $\tilde{\gamma}$ in $\Lambda(n)$, the Maslov index $\widetilde{\mathfrak{m}}(\tilde{\gamma})$ can be calculated as 
\begin{equation}
\widetilde{\mathfrak{m}}(\tilde{\gamma})
= \int_{\tilde{\gamma}} \tilde{\eta}
= \int_{\gamma}d\theta,
\label{eq:maslov index by integration}
\end{equation}
where $\tilde{\eta} = (det_{\mathbb{C}}^{2})^* d\theta$ is the pullback to $\Lambda(n)$ of the canonical $1$-form $d\theta$ on $S^{1}$,
and $\gamma=det_{\mathbb{C}}^{2}(\tilde{\gamma})$ is a loop in $S^{1}$,
see also~\cite{Contreras2003}. 

In the general context of a symplectic manifold $(M,\omega)$ with possibly nontrivial Maslov $S^1$ bundle $\Gamma_J^2$, the integration in Eq.~\eqref{eq:maslov index by integration} can still be defined by treating $\tilde{\eta}$ as the pullback to $\Lambda_J$ of a connection $1$-form on $\Gamma_{J}^{2}$. 
We further note that a global section of $\Gamma_{J}^{2}$ is just an integral manifold of a flat connection (which is an integrable horizontal distribution) on $\Gamma_{J}^{2}$.
To generalize the notion of Maslov indices for the case where $\Gamma_{J}^{2}$
is a nontrivial bundle, we can replace the ``integrable horizontal
distribution'' simply with ``a connection''. 
We thus introduce a nonintegrable version of Maslov indices for smooth loops in $\Lambda_{J}$ with respect to an arbitrary connection $1$-form\footnote{We refer to Appendix~\ref{sec:principal-s1-bundles} for our notation and conventions concerning principal $S^1$ bundles.}
$\alpha=f_\alpha\,\frac{\partial}{\partial\theta}$ on $\Gamma_{J}^{2}$.

\begin{defn}
\label{def:maslov-data} 
Given a connection $1$-form $\alpha = f_\alpha\,\frac{\partial}{\partial\theta}$ on $\Gamma_{J}^{2}$, let $\tilde f_\alpha = (det_J^2)^* f_\alpha$. 
Then for a smooth loop $\tilde{\gamma}:S^1\to\Lambda_{J}$, the \emph{Maslov data of $\tilde{\gamma}$ with respect to $\alpha$} is
\[
\widetilde{\mathfrak{md}}_{\alpha}(\tilde{\gamma}) = \int_{\tilde{\gamma}}\tilde f_\alpha,
\]
and for a smooth loop $\gamma : S^1 \to \Gamma_J^2$ the \emph{Maslov data of $\gamma$ with respect to $\alpha$} is
\begin{equation}
\mathfrak{md}_{\alpha}(\gamma) = \int_{\gamma} f_\alpha.
\label{eq:Maslov-quantity for a path/loop}
\end{equation}
\end{defn}

Clearly, the Maslov data for $\tilde\gamma$ and for $\gamma = det_{J}^{2}\circ\tilde{\gamma}$ are equal:
\[ \widetilde{\mathfrak{md}}_{\alpha}(\tilde{\gamma}) 
= \int_{\tilde{\gamma}} \tilde f_\alpha
= \int_{\tilde{\gamma}} (det_J^2)^* f_\alpha
= \int_{det_J^2 \circ \tilde\gamma} f_\alpha
= \mathfrak{md}_{\alpha}(det_J^2 \circ \tilde\gamma)
= \mathfrak{md}_{\alpha}(\gamma)
, \]
and therefore we will use them interchangeably.

\begin{rem}
We avoid using the term \emph{index} to emphasize that the Maslov
data is not necessarily a topological quantity, even when $\tilde{\gamma}$
is a loop. 
\end{rem}

\begin{rem}
When the horizontal distribution $\mathcal H = \ker\alpha$ associated to the connection $1$-form $\alpha$ has an integral manifold $\mathcal{S}_{0}$ which is an one-sheet covering over $M$, Definition~\ref{def:maslov-data} gives the Maslov indices in the usual sense with respect to the Maslov cycle $\mathfrak{M}_{0}=(det_{J}^{2})^{-1}(\mathcal{S}_{0})$ (or to the section $\mathfrak{s}$). 
Note that this is the case in \cite{Contreras2003}
where a dynamical system is considered in the vicinity of a connected submanifold $\Sigma$ in $M$ on which (the restriction of) the bundle of Lagrangian planes $\Lambda_{pl}$ can be trivialized as $\Sigma\times\Lambda(n)$.
Also note that in general, $\ker\alpha$ being flat does not necessary mean that the integral manifolds are single-sheeted coverings. 
A discussion of nontrivial bundles $\Gamma_{J}^{2}$ with flat connections is beyond the scope of this work, and we do not go into any further details. 
\end{rem}

This work explores the interplay between the dynamics of group actions on the Maslov $S^1$ bundles $\Gamma_J$ and $\Gamma_J^2$, and the geometry of these bundles.
Moreover, it investigates the notion of Maslov data, Definition~\ref{def:maslov-data}, as a nonintegrable version of ordinary Maslov indices. 
We give now a detailed summary and outline of this paper, while introducing the main results.

Section~\ref{sec:basic-properties-maslov-indices} discusses basic properties of Maslov indices. It shows that if the Maslov $S^1$ bundle $\Gamma_J^2$ is trivial, then the Maslov index for a loop $\gamma$ in $M$ is an even number. Moreover, it shows that if $M$ is simply connected and $\Gamma_J^2$ is trivial, the Maslov index does not depend on the choice of section $\mathfrak s^2 : M \to \Gamma_J^2$.

Given the flow $\varphi:\mathbb{R}\times M\rightarrow M$ of an arbitrary symplectic vector field $X$, each tangent map $\varphi_{*}^{t}$ maps a Lagrangian space to another.
Namely, $\varphi_{*}$ defines a flow on $\Lambda_{pl} \cong \Lambda_J$ which covers $\varphi$. 
Since each trajectory $\tilde{\gamma}$ of $\varphi_{*}$ is mapped to a path $\gamma=det_{J}^{2}\circ\tilde{\gamma}$ in $\Gamma_{J}^{2}$, it is natural to expect that, under some conditions, the composition
$
det_{J}^{2}\circ\varphi_{*} : \mathbb{R}\times\Lambda_J \to \Gamma_{J}^{2}
$
factors through $\mathbb{R}\times\Gamma_{J}^{2}$ and induces a flow $\varphi_{\Gamma^{2}}$ on $\Gamma_{J}^{2}$.
Section~\ref{sec:compact-group-actions} shows that such a factorization indeed exists when the vector field $X$ is the infinitesimal generator of a symplectic action by a compact Lie group $G$. 
In particular, it shows that the symplectic action of a compact Lie group $G$ on $(M,\omega)$ can be lifted to the bundles of unitary frames $Fr_J^u$ and to the Maslov $S^1$ bundles $\Gamma_J$ and $\Gamma_J^2$.
The lifted actions commute with the inherent $S^1$ action on each of the Maslov $S^1$ bundles.

An example of a symplectic manifold with nontrivial Maslov $S^1$ bundle is provided by $M=S^2$ with the standard symplectic form $\omega_{S^2}$. After constructing the Maslov $S^1$ bundles $\Gamma_J$ and $\Gamma_J^2$,
Section~\ref{sec:example-S2} considers the symplectic action of $\mathbb{SO}(3)$ on $S^2$ from two different points of view.
First, note that $S^2$ is a symplectic homogeneous space since the symplectic $\mathbb{SO}(3)$ action is transitive. 
The analysis of the $\mathbb{SO}(3)$ action on $S^2$ shows that the lifted $\mathbb{SO}(3)$ actions on the Maslov $S^1$ bundles $\Gamma_J$ and $\Gamma_J^2$ are also transitive, see Proposition \ref{prop:so3-transitive-maslov-bundle}, that is, $\Gamma_J$ and $\Gamma_J^2$ are also homogeneous spaces.
Second, note that $H_{dR}^{1}(S^{2})=0$, hence any symplectic vector field on $S^{2}$ is Hamiltonian.
Proposition~\ref{prop:so3-symplectic-action} asserts that the Hamiltonian for the infinitesimal generators of the $\mathbb{SO}(3)$ action on $M$ can be expressed in terms of an $\mathbb{SO}(3)$-invariant connection $1$-form on the Maslov $S^1$ bundle $\Gamma_J$. 

The properties of the $\mathbb{SO}(3)$ action on $S^2$ discussed in Section~\ref{sec:example-S2} are then extended to more general settings in Section~\ref{sec:maslov-bundles-homogeneous} and Section~\ref{sec:monotone-symplectic}. 
Section~\ref{sec:maslov-bundles-homogeneous} establishes the following generalization of Proposition \ref{prop:so3-transitive-maslov-bundle}.

\begin{thm}
\label{thm:homogeneous-G-space-(intro)}
Let $G$ be a compact Lie group acting transitively and symplectically on $M$, that is, $M$ is homogeneous $G$-space. 
If the Chern class $c_\Gamma$ of the Maslov $S^1$ bundle $\Gamma_J$ is non-zero, then $\Gamma_J$ is also a homogeneous $G$-space.
\end{thm}

Moreover, Section~\ref{sec:monotone-symplectic} generalizes Proposition \ref{prop:so3-symplectic-action} to the case where $(M,\omega)$ is a symplectic manifold such that $[\omega] = r \, c_\Gamma$ for some $r \in \mathbb R$. 
In particular, suppose that $\Phi: G \times M \to M$ is the symplectic action of a compact Lie group $G$ on $M$.
Denote by $\Phi_\Gamma : G \times \Gamma_J \to \Gamma_J$ the lift of $\Phi$ on $\Gamma_J$. 
For a fixed $v \in \mathfrak{g}$, where $\mathfrak{g}$ is the Lie algebra of $G$, denote by $X_v$ the infinitesimal generator of $\Phi$ associated to $v$, and by $\mathcal{X}_{v}$ the corresponding infinitesimal generator of $\Phi_\Gamma$. 
We have the following result.

\begin{thm}\cite{Ono1988, Lupton1995}
\label{thm:hamiltonian-G-action}
Let $(M,\omega)$ be a symplectic manifold, and denote by $c_{\Gamma}\in H_{dR}^{2}(M)$ the first real Chern class of the Maslov $S^1$ bundle $\pi_\Gamma : \Gamma_J \to M$.
If $[\omega]=r \, c_{\Gamma}$ for some real number $r$, then any symplectic action $\Phi$ on $M$ by a compact Lie group $G$ is Hamiltonian. 
More specifically, for any $v\in\mathfrak{g}$, there exists a Hamiltonian $H_v$ for the symplectic vector field $X_v$, satisfying $\pi_{\Gamma}^*(H_v) = \bar{\eta}(\mathcal{X}_{v})$ for a $G$-invariant $1$-form $\bar{\eta}$ with $d\bar{\eta}=-\pi_\Gamma^* \omega$. 
\end{thm}

\begin{rem}
The conclusion in Theorem \ref{thm:hamiltonian-G-action} about symplectic actions being Hamiltonian has been obtained in \cite{Ono1988} and \cite{Lupton1995} for circle actions on compact monotone manifolds, and as is mentioned in \cite{Cho2017}, it sufficiently implies the conclusion for symplectic actions by compact Lie groups. 
While we work with different bundle structures and we do not assume that the symplectic manifold is compact, our proof of Theorem \ref{thm:hamiltonian-G-action} is in the same spirit with the corresponding proof in \cite{Ono1988}. 
Nevertheless, our demonstration provides a clearer connection between the geometry related to the Maslov index and the Hamiltonian dynamics. 
Note that \cite{Ono1988, Ono1992, Lupton1995, Cho2017} focus on the obstructions to the existence of Hamiltonian group actions, for which the compactness of the manifold $M$ plays a crucial role.
We further discuss this in Section~\ref{sec:cotangent-bundles}.
\end{rem}

Section~\ref{sec:maslov-data-s1} considers symplectic circle actions focusing on their Maslov data. 
When $G=S^1$, the orbits of the lifted action $\Phi_{\Gamma}$ on $\Gamma^2_{J}$ covering the same orbit of $\Phi$ on $M$ have the same Maslov data. 
In other words, for $w,w'\in\pi_\Gamma^{-1}(p)$ with $\gamma_{w}(z)=\Phi_{\Gamma}^{z}(w)$ and $\gamma_{w'}(z)=\Phi_{\Gamma}^{z}(w')$ for all $z \in S^1$, we have $\mathfrak{m}_{\alpha}(\gamma_{w}) = \mathfrak{m}_{\alpha}(\gamma_{w'})$ for any connection $1$-form $\alpha$ on $\Gamma^2_J$. 
Therefore, this defines a smooth function $Q_\alpha$ on $M$ with $Q_{\alpha}(p) = \mathfrak{m}_{\beta}(\gamma_{w})$, which we call the Maslov data of $\Phi$ with respect to $\alpha$.
Although with a different connection $\alpha$ the function $Q_\alpha$ may be different, its values at the fixed points of the action $\Phi$ turn out to be independent of $\alpha$. 
We call these values the \emph{local Maslov indices} at the corresponding fixed points.
Then, Proposition \ref{prop:maslov-index-phi} states that when the global Maslov indices are defined, they should have the same value as the local Maslov indices, which is independent of the choice of the global section of $\Gamma_J^2$. 

In Section~\ref{sec:int-ham-sys} we turn our attention to integrable Hamiltonian systems of $2$ degrees of freedom. 
First, Proposition \ref{prop:action-angle coordinates and local Maslov indices} discusses conditions under which an $S^1$ symmetry of an integrable Hamiltonian system can be extended to a $\mathbb T^2$ symmetry. 
Then, Section~\ref{sec:pinched-tori} shows that the restriction of $\Gamma_J^2$ over a Lagrangian pinched torus, associated to focus-focus singularities and Hamiltonian monodromy \cite{Duistermaat1980}, is a trivial bundle. 

Section~\ref{sec:S2xS2} then considers the construction of the Maslov $S^1$ bundle $\Gamma_{S^2 \times S^2}$ over the symplectic manifold $S^2 \times S^2$ which appears in applications from physics such as coupled spin problems and the Kepler problem, and it computes the local Maslov indices at the fixed points of the $S^1$ action that rotates the two spheres.

Finally, we provide a brief summary of the results in Section~\ref{sec:conclusions} and collect auxiliary definitions and results in a series of appendices.
In particular, Appendix~\ref{sec:notation} recalls basic facts and fixes notation for principal $S^1$ bundles and unitary frames, while Appendix~\ref{sec:topology-lambda} discusses the topology of $\Lambda_{pl} \cong \Lambda_J$, and Appendix~\ref{sec:structure-unitary-frame-bundle} discusses the structure of the unitary frame bundle $Fr^u_J$.

\section{Basic properties of Maslov indices}
\label{sec:Basics of Maslov-Index}
\label{sec:basic-properties-maslov-indices}

In this section we establish some basic properties of Maslov indices over Lagrangian subbundles and we show that, if the manifold $M$ is simply connected, then the Maslov index $\mathfrak{m}_\mathfrak{s}$ does not depend on the choice of the section $\mathfrak{s} : M \to \Gamma_J^2$.

Let $\mathcal{N}$ be a submanifold in $M$ and $\mathfrak{L}_{\mathcal{N}}$ be a Lagrangian subbundle of $TM$ over $\mathcal{N}$. 
That is, $\mathfrak{L}_{\mathcal{N}}\big|_{b}$
is a Lagrangian subspace of $T_b M$ for any $b\in\mathcal{N}$. 
Then $\mathfrak{L}_{\mathcal{N}}$ induces a smooth section
\[ \sigma_{\mathcal{N}} : \mathcal{N} \to \Lambda_{J}\big|_{\mathcal{N}} = Fr_{J}^{u}\big/\mathbb{O}(n)\big|_{\mathcal{N}}
: b \mapsto \mathfrak{L}_{\mathcal{N}}\big|_{b} \] 
of the bundle of Lagrangian planes over $\mathcal{N}$.
We have the following result.

\begin{prop}
\label{prop:section-orientable-subbundle}
If $\mathfrak{L}_{\mathcal{N}}$ is an orientable Lagrangian subbundle over a submanifold $\mathcal{N} \subseteq M$, 
then it induces a section $\tilde{\sigma}_{\mathcal{N}} : \mathcal{N} \to Fr_{J}^{u}\big/\mathbb{SO}(n)\big|_{\mathcal{N}}$ of the bundle of oriented Lagrangian planes over $\mathcal{N}$.
\end{prop}

\begin{proof}
Since the bundle $\mathfrak{L}_{\mathcal{N}}$ is orientable, we fix an orientation. 
For any $b\in\mathcal{N}$, there is a neighbourhood $U_{b}$ in $\mathcal{N}$ on which it admits an ordered local frame
$f_{U_{b}}^{*}=(e_{1},...,e_{n})$ of $\mathfrak{L}_{\mathcal{N}}$ with the same orientation. 
With some modification, this can be extended to an ordered unitary frame
\[
f_{U_{b}} = (e_{1},...,e_{n},s_{1},...,s_{n})
\]
for $TM$ on $U_{b}$ with $s_{i}=Je_{i}$. 
Then,
\[
\sigma_{b}:U_{b} \to Fr_{J}^{u} : x \mapsto f_{U_{b}}(x) 
\]
is a local section for $Fr_{J}^{u}$. 
By composing with the quotient map $q_{\mathbb{SO}(n)} :Fr_{J}^{u} \to Fr_{J}^{u}\big/\mathbb{SO}(n)$ we get a local section $\tilde{\sigma}_{\mathcal{N},b}= q_{\mathbb{SO}(n)}\circ\sigma_{b} : U_{b} \to Fr_{J}^{u}\big/\mathbb{SO}(n)$. 

The collection $\{U_{b}\big|b\in\mathcal{N}\}$ constitutes an open cover of $\mathcal{N}$ over each element of which there is a smooth section $\sigma_{\mathcal{N},b}$. 
It remains to check that when $U_{b}\cap U_{b'}\neq\emptyset$, $\tilde{\sigma}_{\mathcal{S},b}$ agrees with $\tilde{\sigma}_{\mathcal{S},b'}$ on $U_{b}\cap U_{b'}$.
Since both $f_{U_{b}}^{*}$ and $f_{U_{b'}}^{*}=(e'_{1},...,e'_{n})$
can be taken as orthonormal frames with respect to the Riemannian metric
$g_{J}$ that fit the same orientation, at each $x\in U_{b}\cap U_{b'}$,
they are related by a matrix $C_x \in \mathbb{SO}(n)$ via
\[
(e'_{1},...,e'_{n}) \big|_x =(e_{1},...,e_{n}) \big|_x \cdot C_x.
\]
For $f_{U_{b'}}=(e'_{1},...,e'_{n},s'_{1},...,s'_{n})$, it holds
\[
s'_{i}=J e'_{i},
\]
and then it is straightforward to check that at each $x \in U_{b}\cap U_{b'}$ we have
\[
(e'_{1},...,e'_{n},s'_{1},...,s'_{n}) \big|_x =(e_{1},...,e_{n},s_{1},...,s_{n}) \big|_x \cdot C_x.
\]
This implies that $\tilde{\sigma}_{\mathcal{N},b'}(x) = \tilde{\sigma}_{\mathcal{N},b}(x)$.
Then we get $\tilde{\sigma}_{\mathcal{N}}$ by piecing the
local sections $\tilde{\sigma}_{\mathcal{N},b}$ together.
\end{proof}

Consider now the case where $\mathcal{N}=M$ and $\mathcal{L}_{\mathcal{N}} = T\mathcal{S}$, with $\mathcal{S}\subset M$ an oriented Lagrangian submanifold of $M$.

\begin{prop}
Assume that $\mathcal{S}$ is an orientable Lagrangian submanifold of $M$ and that the Maslov $S^1$ bundle $\Gamma_{J}$ admits a global section $\mathfrak{s} : M \to \Gamma_J$. Then for any loop $\gamma$ in $\mathcal{S}$, the Maslov index $\mathfrak{m}_{\mathfrak{s^{2}}}(\gamma)$ with respect to the associated global section $\mathfrak{s}^2=q_{\pm1} \circ \mathfrak{s} : M \to \Gamma_J^2$
is an even number; here, $q_{\pm1}:\Gamma_J \to \Gamma_J^2$ is the corresponding quotient map.
\end{prop}

\begin{proof}
The existence of the global section $\mathfrak{s}^2 = q_{\pm1} \circ \mathfrak{s} : M \to \Gamma_J^2$ follows directly from the definitions of $\mathfrak{s} : M \to \Gamma_J$ and $q_{\pm1} : \Gamma_J \to \Gamma_J^2$. We denote by 
$tr_{\mathfrak{s}}:\Gamma_J \to M\times S^1$ and 
$tr_{\mathfrak{s}^2}:\Gamma_J^2 \to M\times S^1$ the trivializations induced by $\mathfrak{s}$ and $\mathfrak{s}^2$ respectively.
Moreover, since $\mathcal{S}$ is orientable, $T\mathcal{S}$ is also orientable, and Proposition \ref{prop:section-orientable-subbundle} implies that there is a global section $\tilde{\sigma}_{\mathcal{S}} : \mathcal{S} \to Fr_J^u \big/ \mathbb{SO}(n)$.

Then the conclusion that the Maslov index $\mathfrak{m}_{\mathfrak{s^{2}}}(\gamma)$ is an even number can be read off the following commutative diagram
\begin{equation}
\begin{CD}
S^{1} 
@> \tilde{\sigma}_{\mathcal{S}}\circ\gamma >> Fr_{J}^{u}\big/\mathbb{SO}(n)
@> det_{J} >> \Gamma_{J}
@> tr_{\mathfrak{s}} >> M\times S^{1} 
@> pr_{S^{1}} >> S^{1} 
\\
@| 
@Vq_{\pm1}VV 
@Vq_{\pm1}VV 
@Vq_{\pm1}VV 
@V{\scriptstyle \mathrm{square}}VV
\\
S^{1}
@> \sigma_{\mathcal{S}}\circ\gamma >> \Lambda_{J}
@> det_{J}^{2} >> \Gamma_{J}^{2}
@> tr_{\mathfrak{s}^{2}} >> M\times S^{1}
@> pr_{S^{1}} >> S^{1}
\end{CD}
\label{Diagram:oriented Maslov indices and Maslov indices}
\end{equation}
which shows that the degree of the map at the bottom row is twice the degree of the map at the top row.
\end{proof}

\begin{rem}
A Lagrangian subbundle $\mathfrak{L}$ of $TM$ over $M$ induces a global section $\sigma:M\rightarrow\Lambda_{J}$ by sending each $b\in M$ to the fiber $\mathfrak{L}_{b}$ of $\mathfrak{L}$ above it. 
Then $\mathfrak{s}^{2}:=det_{J}^{2}\circ\sigma$ is a global section to the bundle $\Gamma_{J}^{2}$. 
Furthermore, if $\mathfrak{L}$ is orientable, then it also induces a trivialization
for $\Gamma_{J}$.
\end{rem}

Generally speaking, a loop $\tilde{\gamma} : S^1 \to \Lambda_J$ of Lagrangian planes may have different Maslov indices with respect to different global sections of $\Gamma_{J}^{2}$.
However, when the manifold $M$ is simply connected we obtain the following result.

\begin{prop}
\label{prop:maslov-index-simply-connected}
Given a simply connected symplectic manifold $M$, the Maslov index $\mathfrak{m}_{\mathfrak{s}^2}(\tilde \gamma)$ of a loop $\tilde\gamma: S^1 \to \Lambda_J$ does not depend on the choice of a global section $\mathfrak{s}^2$ of $\Gamma_J^2$.
\end{prop}

\begin{proof}
Suppose that $\mathfrak{s}^2$ and $\hat{\mathfrak{s}}^2$ are sections of $\Gamma_J^2$.
In terms of the trivialization $\Gamma_J^2\xrightarrow{tr_{\hat{\mathfrak{s}}^2}}M\times S^{1}$ with respect to $\hat{\mathfrak{s}}^2$, the section $\mathfrak{s}^2:M\rightarrow M\times S^{1}$ takes the form $\mathfrak{s}^2(p)=(p,\theta_{p})$ with $\theta:p\mapsto\theta_{p}$ being a smooth map from $M$ to $S^{1}$. 
A loop $\gamma : S^1 \to \Gamma_J^2 : z \mapsto \gamma(z)$ takes the form $z\mapsto\big(\lambda(z),\tau'(z)\big)$ in the trivialization $\Gamma_J^{2}\xrightarrow{tr_{\hat{\mathfrak{s}}^2}}M\times S^{1}$, and the form $z\mapsto\big(\lambda(z),\tau(z)\big)$ in $\Gamma_J^{2}\xrightarrow{tr_{\mathfrak{s}^2}}M\times S^{1}$, with maps $\lambda:S^{1}\rightarrow M$, $\tau',\tau:S^{1}\rightarrow S^{1}$.
The corresponding Maslov indices $\mathfrak{m}_{\hat{\mathfrak{s}}^2}(\gamma)$ and $\mathfrak{m}_{\mathfrak{s}^2}(\gamma)$ are then the degrees of
the maps $\tau'$ and $\tau$, respectively. 
It holds that $\tau'(z)=\tau(z)\cdot\theta_{\lambda(z)}$.
Since $M$ is simply connected, the mapping $z\mapsto\theta_{\lambda(z)}$ which is subject to the factorization $S^{1}\xrightarrow{\lambda}M\xrightarrow{\theta}S^{1}$ has degree $0$, and hence $\tau$ and $\tau'$ have the same degree.
\end{proof}

\section{Compact group actions on Maslov \texorpdfstring{$S^{1}$}{S1} bundles}
\label{sec:compact-group-actions}

In this section, we show that the symplectic action $\Phi$ of a compact Lie group $G$ on a symplectic manifold $(M,\omega)$ can be lifted to a $G$-action on $\Gamma_J$ (respectively, $\Gamma_J^2$) which covers $\Phi$ and commutes with the inherent $S^{1}$ action on $\Gamma_{J}$ (respectively, $\Gamma_J^2$).
In particular, let
\[
\Phi : G \times M \to M
\]
be a symplectic left action on $M$ by a compact Lie group $G$.
Namely, $\Phi$ satisfies that $\Phi^{h'h} = \Phi^{h'} \circ \Phi^{h}$ and $(\Phi^{h})^*\omega = \omega$ for all $h, h' \in G$. 

\begin{defn}
\label{def:group action on the bundle}
A \emph{group action by $G$} on a principal $S^1$ bundle $P\rightarrow B$ is a $G$-action on the manifold $P$ which commutes with the inherent $S^{1}$-action of the principal bundle.
\end{defn}

The symplectic $G$-action $\Phi$ on $M$ can be lifted to a $G$-action on the bundle $Fr^u_J \to M$ by resorting to a $G$-equivariant almost complex structure $\bar J$.
The latter can be constructed by first averaging an arbitrary Riemannian metric $g$ over $G$ to obtain a $G$-invariant metric $\bar g$, and then defining $\bar J$ through the polar decomposition of the $TM$ isomorphism $\mathcal A$ satisfying $\omega(u,v) = \bar{g}(\mathcal{A} u,v)$, see \cite[Lemma 5.5.6]{McDuff2017}.
Then, $g_{\bar J}(\cdot,\cdot)=\omega(\bar J \cdot, \cdot)$ is a Riemannian metric compatible with $\omega$ and $\bar J$, while $\Phi^h_* \circ \bar J = \bar J \circ \Phi^h_*$ and $(\Phi^h)^* g_{\bar J} = g_{\bar J}$.

The symplectic action $\Phi$ of $G$ on $(M,\omega)$ lifts naturally to the following group action by $G$ on the principal bundle $Fr^u_{\bar J} \to M$:
\[
\Phi_{\#}:G\times Fr_{\bar{J}}^{u}\rightarrow Fr_{\bar{J}}^{u}
: h, (u_1,\dots,v_n) \mapsto 
\Phi_{\#}^{h}(u_1,\dots,v_n)
= \big(\Phi^h_*(u_1), \dots, \Phi^h_*(v_n)\big),
\]
where $(u_1,\dots,u_n,v_1,\dots,v_n)$ denotes a unitary frame with respect to $\bar J$, see Appendix~\ref{sec:unitary-frames}, and it can be shown that $\big(\Phi^h_*(u_1), \dots, \Phi^h_*(v_n)\big)$ is also unitary with respect to $\bar J$.
To show that $\Phi_\#$ commutes with the inherent $\mathbb{U}(n)$ action of the principal bundle $Fr^u_{\bar J} \to M$, it can be checked that
\begin{equation}
\Phi^h_\#((u_1,\dots,v_n) \cdot C)
= \big(\Phi^h_*(u_1), \dots, \Phi^h_*(v_n)\big) \cdot C
= \Phi_{\#}^{h}(u_1,\dots,v_n) \cdot C,
\label{eq:phi-sharp-equivariance}
\end{equation}
for any $C \in \mathbb{U}(n)$.
Hence, $\Phi_{\#}$ is a $G$-action on $Fr_{\bar{J}}^{u} \to M$.

Due to Eq.~\eqref{eq:phi-sharp-equivariance}, $\Phi_{\#}$ induces a smooth $G$-action $\Phi_\Gamma$ on $\Gamma_{\bar{J}}=Fr_{\bar{J}}^{u}\big/\mathbb{SU}(n)$.
More specifically, let $[u_1,\dots,v_n]$ denote the equivalence class of the unitary frame $(u_1,\dots,v_n)$ under the right $\mathbb{SU}(n)$ action on $Fr^u_{\bar J}$.
Then, the action $\Phi_\Gamma$ is defined as
\begin{equation}
\Phi_\Gamma 
: G\times \Gamma_{\bar J} \to \Gamma_{\bar J}
: h, [u_1,\dots,v_n] \mapsto 
[\Phi_{\#}^{h}(u_1,\dots,v_n)].
\end{equation}
Moreover, $\Phi_{\Gamma}$ on $\Gamma_{\bar{J}}$ commutes with the $S^{1}$ action. 
To see this, note that the $S^1 \cong \mathbb{U}(n) \big/ \mathbb{SU}(n)$ action on $\Gamma_{\bar J}$ is defined by
\[
[u_1,\dots,v_n] \cdot e^{2\pi i\theta}
:= [u_1,\dots,v_n] \cdot [C] 
= [(u_1,\dots,v_n) \cdot C],
\]
where $e^{2\pi i\theta} = det_{\mathbb{C}}(C)$.
Then, we check that
\[
\begin{aligned}
\Phi_{\Gamma}^{h} \big( [u_1,\dots,v_n]\cdot e^{2 \pi i\theta} \big) 
& = \Phi_{\Gamma}^{h} \big( [u_1,\dots,v_n] \cdot [C] \big)
= \Phi_{\Gamma}^{h} \big( [(u_1,\dots,v_n) \cdot C] \big)
= \big[ \Phi_\#^h \big( (u_1,\dots,v_n) \cdot C \big) \big] \\
& = \big[ \Phi_\#^h (u_1,\dots,v_n) \cdot C \big]
= [ \Phi_\#^h (u_1,\dots,v_n) ] \cdot [C] 
= \Phi_\Gamma^h[u_1,\dots,v_n] \cdot [C] \\
& = \Phi_\Gamma^h[u_1,\dots,v_n] \cdot e^{2 \pi i\theta}.
\end{aligned}
\]

The same argument can be used to define a group action $\Phi_{\Gamma^2}$ by $G$ on the bundle $\Gamma_{\bar J}^2 = Fr^u_{\bar J} \big/ \mathbb{EU}(n) \to M$.
Since the unitary structure on $M$ is unique up to isomorphism, the
discussion in this subsection amounts to the following proposition.

\begin{prop}
\label{prop:G-action-lift}
Suppose that $\Phi$ is a symplectic group action on $(M,\omega)$ by a compact Lie group $G$ and $J$ is a compatible almost complex structure.
Then $\Phi$ lifts to a group action $\Phi_\#$ by $G$ on the principal bundle $Fr_{J}^{u} \to M$ that covers $\Phi$,
and it induces group actions $\Phi_\Gamma$ and $\Phi_{\Gamma^2}$ by $G$ on the principal bundles $\Gamma_J \to M$ and $\Gamma_J^2 \to M$, respectively.
The actions $\Phi_\Gamma$ and $\Phi_{\Gamma^2}$ cover $\Phi$ and are covered by $\Phi_{\#}$.
\end{prop}

For each $p\in M$ and $w\in\Gamma_{J}$, denote by $G_{p}$ the isotropy group of the action $\Phi$ at $p$, and by $G_{w}$ the isotropy group of the lifted action $\Phi_{\Gamma}$ at $w$.
Let $w$ be a point in $\pi_{\Gamma_{J}}^{-1}(p)$. 
Then $G_w \subseteq G_p$.
Since the action $\Phi_{\Gamma}$ on $\Gamma_{J}$ covers the action $\Phi$ on $M$, $G_p$ acts on the fiber $\pi_{\Gamma_{J}}^{-1}(p)$. 
That is, $h\cdot w\in\pi_{\Gamma_{J}}^{-1}(p)$ for any $h\in G_p$ and any $w\in\pi_{\Gamma_{J}}^{-1}(p)$.
It turns out that $G_{w'}=G_w$ for any $w'\in\pi_{\Gamma_{J}}^{-1}(p)$, and that the $G_p$-action on $\pi_{\Gamma_{J}}^{-1}(p)$ induces a homomorphism $\phi_p$ from $G_p$ to $S^{1}$ with kernel $G_w$. 

\begin{prop}
\label{prop:phi_p}
For each $p \in M$ there exists a Lie group homomorphism $\phi_p: G_p \to S^1$ with kernel $G_w$  such that, for all $w\in\pi_{\Gamma_{J}}^{-1}(p)$ and $h\in G_{p}$, it holds
\[
h\cdot w=w\cdot\phi_{p}(h).
\]
Moreover, the family $\{\mathcal{\phi}_{p}\}$ is $G$-related in
the sense that, for all $p\in M$ and $g\in G$, it holds
\begin{equation}
\phi_{g\cdot p}\circ Ad_{g}=\phi_{p}.
\label{eq:adjoint-action-phi-p}
\end{equation}
\end{prop}

\begin{proof}
Fix points $p \in M$ and $w \in \pi_{\Gamma_J}^{-1}(p)$. 
Since the inherent $S^1$ action on $\Gamma_J$ acts transitively and freely on $\pi_{\Gamma_J}^{-1}(p)$, for each $h\in G_p$ there exists a unique $z_h \in S^1$ such that $h \cdot w = w \cdot z_{h}$. Then 
\[
w\cdot z_{h'h}=(h'h)\cdot w = w \cdot z_{h} \cdot z_{h'} = w \cdot(z_{h'}z_{h}),
\]
and hence $z_{h'h}=z_{h'}z_{h}.$ Define $\phi_p$ as $\phi_p(h)=z_h$.
Then $\phi_p$ is a homomorphism, and $h\cdot w = w$ if and only if $\phi_p(h)=1_{S^1}$. 
Hence $\ker\phi_p=G_w$.
Note that $\pi_{\Gamma_{J}}^{-1}(p)$ is diffeomorphic to $S^{1}$
via $\bar{\theta}_{w}:w\cdot z\mapsto z$. Since $\phi_{p}(h)=\bar{\theta}_{w}\circ\Phi_{\Gamma}^{h}(w)$,
$\phi_{p}$ is a smooth map, and therefore it is a Lie group homomorphism.

The assignment of $z_{h}$ to $h$ is independent of the choice of $w$. 
Namely, the identity $h\cdot w'=w'\cdot z_{h}$ holds for all $w'\in\pi_{\Gamma_{J}}^{-1}(p)$. 
To see this, first note that there exists $z\in S^{1}$ such that $w\cdot z=w'$, and then
\[
h\cdot w'=h\cdot(w\cdot z)=(h\cdot w)\cdot z=w\cdot z_{h}\cdot z=w'\cdot z_{h}.
\]

It remains to show that $G$ acts on the family $\{\mathcal{\phi}_{p}\}$ according to Eq.~\eqref{eq:adjoint-action-phi-p}.
Let $w$ be a point on the fiber $\pi_{\Gamma}^{-1}(p)$. Then $w'=g\cdot w\in\pi_{\Gamma}^{-1}(g\cdot p)$.
For any $h\in G_{p}$,
\[
ghg^{-1}\cdot w'
= ghw\\
= g\cdot\big(w\cdot\phi_{p}(h)\big)\\
= (g\cdot w)\cdot\phi_{p}(h)\\
= w'\cdot\phi_{p}(h),
\]
which implies $\phi_{g\cdot p}(ghg^{-1})=\phi_{p}(h)$ and concludes the proof.
\end{proof}
Due to Proposition~\ref{prop:phi_p}, $\operatorname{im}\phi_p = \phi_p(G_p) \cong G_p / G_w$ is a compact subgroup of $S^1$ and therefore $\operatorname{im}\phi_p$ is either a finite cyclic group $\mathbb Z_k = \{e^{2\pi i l/k} \,|\, l=0,\dots,k-1\}$, or $S^{1}$ itself.
Moreover, due to Eq.~\eqref{eq:adjoint-action-phi-p}, for any $p$, $p'$ on the same $G$-orbit, we have $\operatorname{im}\phi_{p}=\operatorname{im}\phi_{p'}$.
In particular, if $G$ acts transitively on $M$, then $\operatorname{im}\phi_{p}=\operatorname{im}\phi_{p'}$ for all $p, p' \in M$.

\section{An example with nontrivial Maslov bundle: \texorpdfstring{$S^2$}{S2}}
\label{subsec:An-example: S^2}
\label{sec:example-S2}

In this section we consider the case of the symplectic manifold $M=S^2$ which may be the simplest example with non-trivial Maslov $S^1$ bundles $\Gamma_J$ and $\Gamma_J^2$.
After describing $\Gamma_J$, we consider the action of $\mathbb{SO}(3)$ to $S^2$ and we show that it lifts to a transitive Hamiltonian action on $\Gamma_J$. 
These results are special instances of the more general results obtained in the subsequent two sections. 

Consider $S^2$ as the embedded submanifold $\{(x,y,z): x^2+y^2+z^2=1\}$ of $\mathbb{R}^3$.
At each point $p \in S^{2}$, the tangent space $T_{p}S^{2}$ is a
subspace of $T_{p}\mathbb{R}^{3}$. 
Let $n$ be the restriction of the vector field $x\frac{\partial}{\partial x}+y\frac{\partial}{\partial y}+z\frac{\partial}{\partial z}$
on $S^{2}$. 
With $T_{p}\mathbb{R}^{3}$ being identified with $\mathbb{R}^{3}$, it holds $n_p=p$. 
Then,
\begin{equation}
\label{eq:omega_S2}
\omega_{S^{2}} = \iota_{n} dx \wedge dy \wedge dz
\end{equation}
is a symplectic structure on $S^{2}$. 
Viewing $u, v \in T_{p}S^{2}$ as vectors in $\mathbb{R}^{3}$, it holds
\[
\omega_{S^{2}}(u,v) = p \cdot (u\times v).
\]
With respect to $\omega_{S^{2}}$, 
\[
J_{S^{2}}(u):=-p\times u
\]
defines a compatible almost complex structure, and 
\[
g_{S^{2}}(u,v):=\omega_{S^{2}}(J_{S^{2}}u,v) = u \cdot v
\]
is the restriction to $S^2$ of the standard Riemannian metric on $\mathbb{R}^3$. 

Denote by $\Gamma_{S^2} := \Gamma_{J}$ and $Fr_{S^2}^{u} := Fr_{J}^{u}$, respectively, the Maslov $S^{1}$ bundle and the unitary frame bundle of $(S^{2},\omega_{S^{2}})$.
Over each $p\in S^{2}$, each unitary frame of $T_p S^2$ can be determined uniquely and distinctly by a unit vector $u\in\mathbb{R}^{3}$ tangent to $S^{2}$ at $p$.
The unitary frame associated to the vector $u$ is denoted by\footnote{For our notation and conventions concerning unitary frames we refer to Appendix~\ref{sec:unitary-frames}.}
\[
\big(u, -J_{S^{2}}(u) ;\,p\big)
= (u, p\times u ;\, p).
\]
The matrix $[u, p\times u, p]$ with $u$ , $p\times u$, and $p$ being the columns, is an element of $\mathbb{SO}(3)$. 
More precisely, the map
\[
\operatorname{mat}:Fr_{S^{2}}^{u} \to \mathbb{SO}(3) : \big(u,-J_{S^{2}}(u);p\big)\mapsto[u,p\times u,p]
\]
is a diffeomorphism from $Fr_{S^{2}}^{u}$ to $\mathbb{SO}(3)$. 
Since $\mathbb{SU}(1)=\{1\}$, $\Gamma_{S^{2}}=Fr_{S^{2}}^{u} \cong \mathbb{SO}(3)$.

The inherent $S^1$ action on the principal bundle $\Gamma_{S^2}$ is given by
\begin{equation}
\label{eq:gamma-s2-inherent-s1-action}
(u, v) \cdot e^{i\theta} 
= (u,v) \cdot \begin{bmatrix}
\cos\theta & -\sin\theta \\ \sin\theta & \cos\theta
\end{bmatrix}
= (u\cos\theta + v \sin\theta, -u \sin\theta + v\cos\theta),
\end{equation}
where $(u,v)$ is a unitary frame at $p \in S^2$ with $v = p \times u$.

Since $\mathbb{SO}(3)$ acts symplectically on $(S^{2},\omega_{S^{2}})$, Proposition~\ref{prop:G-action-lift} implies that there is an $\mathbb{SO}(3)$-action on the bundle $\pi_{\Gamma_{S^{2}}} : \Gamma_{S^{2}} \to S^{2}$ which covers the action on $S^{2}$. 
In fact, this action is exactly given by 
\[
\mathbb{SO}(3)\times\Gamma_{S^{2}}
\xrightarrow{(\operatorname{id},\operatorname{mat})}
\mathbb{SO}(3) \times \mathbb{SO}(3)
\xrightarrow{\operatorname{mult}} \mathbb{SO}(3)
\xrightarrow{\operatorname{mat}^{-1}} \Gamma_{S^{2}}
\]
with
\[
\big(A,[u,p\times u,p]\big)\mapsto
A [u,p\times u,p] = [Au,Ap\times Au,Ap].
\]
Here $\mathbb{SO}(3)\times\mathbb{SO}(3)\xrightarrow{\operatorname{mult}}\mathbb{SO}(3)$
is the group multiplication $(A,B)\mapsto A B$, and hence
the lifted group action is recognized as the left action of $\mathbb{SO}(3)$
on itself (identified as $\Gamma_{S^{2}}$). Since the left action of $\mathbb{SO}(3)$ on itself is transitive, we have the following proposition.

\begin{prop}
\label{prop:SO(3)-action-lift}
\label{prop:so3-transitive-maslov-bundle}
The lifted $\mathbb{SO}(3)$ action is transitive on $\Gamma_{S^{2}}$.
\end{prop}

\begin{rem}
An alternative way of viewing Proposition~\ref{prop:so3-transitive-maslov-bundle} is the following. 
Since $\mathbb{SO}(3)$ acts transitively on the manifold $M=S^2$, the latter is a homogeneous space with $S^2 \cong \mathbb{SO}(3) \big/ \mathbb{SO}(2)$. 
Then Proposition~\ref{prop:so3-transitive-maslov-bundle} asserts that the Maslov $S^1$ bundle  $\Gamma_{S^2}$ is also a homogeneous space with $\Gamma_{S^2} \cong \mathbb{SO}(3) \big/ \{1\}$. 
In Section~\ref{sec:maslov-bundles-homogeneous} we show that manifolds $M$ that are homogeneous $G$-spaces give rise to Maslov $S^1$ bundles that are also homogeneous $G$-spaces.
\end{rem}

We now turn our attention to showing that the symplectic action of $\mathbb{SO}(3)$ on $S^2$ is Hamiltonian.
Let $\alpha=f_{\alpha}\,\frac{\partial}{\partial\theta}$ be a connection $1$-form on $\Gamma_{S^{2}}$. 
Since $S^{2}$ is simply connected and the principal bundle $\pi_{\Gamma_{S^{2}}}:\Gamma_{S^{2}}\rightarrow S^{2}$ is not trivial, there is no integrable connection on $\Gamma_{S^{2}}$ and hence the curvature form $\Omega_{\alpha}$, satisfying $\pi_{\Gamma_{S^{2}}}^* \Omega_\alpha = df_\alpha$, does not vanish, see Appendix~\ref{sec:principal-s1-bundles}. 

An arbitrary connection $1$-form $\alpha$ can be averaged over $\mathbb{SO}(3)$ to produce an $\mathbb{SO}(3)$-invariant connection $1$-form $\bar\alpha$.
In particular, starting with an arbitrary connection $1$-form $\alpha = f_\alpha \, \frac{\partial}{\partial\theta}$, for each element $v$ of the tangent bundle $T\Gamma_{S^2}$ define
\[
\bar{f}(v):=\int_{\mathbb{SO}(3)}f_{\alpha}(h_{*}v)d\mu(h)
\]
with $dh$ being a right invariant probability measure on $\mathbb{SO}(3)$.
Then $\bar{f}$ is both $S^{1}$-invariant and $\mathbb{SO}(3)$-invariant, and $\bar{f}(\frac{\partial}{\partial\theta})\equiv1$. 
Therefore, $\bar{\alpha} = f_{\bar\alpha} \,\frac{\partial}{\partial\theta}$ with $f_{\bar\alpha} := \bar f$ defines a connection $1$-form which is invariant under the $\mathbb{SO}(3)$-action.
Then $df_{\bar\alpha}$ is also $\mathbb{SO}(3)$-invariant, and so is $\Omega_{\bar{\alpha}}$ with $df_{\bar\alpha} = \pi_{\Gamma_{S^{2}}}^{*} \Omega_{\bar{\alpha}}$.

We have the following result.
\begin{lem}
\label{lem:symp-potential on S^2}
If $\bar{\alpha}$ is an $\mathbb{SO}(3)$-invariant connection $1$-form on $\Gamma_{S^{2}}$, then $\Omega_{\bar\alpha} = \omega_{S^2}$. 
\end{lem}

\begin{proof}
Since $S^2$ is $2$ dimensional, there exists a real-valued function $g$ on $S^2$ such that $\Omega_{\bar\alpha} = g \,\omega_{S^2}$. 
Moreover, since $\Omega_{\bar\alpha}$ and $\omega_{S^2}$ are $\mathbb{SO}(3)$-invariant, the function $g$ must also be $\mathbb{SO}(3)$-invariant, and thus constant: $g(p) = r$ for all $p \in S^2$.
Therefore, $\Omega_{\bar\alpha} = r\,\omega_{S^2}$ and thus $[\Omega_{\bar{\alpha}}] = r\,[\omega_{S^2}]$.
Since the first Chern number of the bundle $\Gamma_{S^2} \cong \mathbb{SO}(3) \to S^2$ equals $2$, and $\frac{1}{2\pi} \int_{S^2} \omega_{S^2} = 2$, we conclude that $r = 1$.
\end{proof}

For $v \in \mathfrak{so}(3)$, the corresponding infinitesimal generator for the $\mathbb{SO}(3)$ action on $S^2$ is denoted by $X_v$, while the infinitesimal generator for the $\mathbb{SO}(3)$ action on $\Gamma_{S^2}$ is denoted by $\mathcal{X}_v$.
The following result, which will be generalized in Section~\ref{sec:monotone-symplectic},
asserts that the $\mathbb{SO}(3)$ action on $S^2$ is Hamiltonian and gives an expression for the Hamiltonian function $H_v$ associated to $v \in \mathfrak{so}(3)$ in terms of an $\mathbb{SO}(3)$-invariant connection $1$-form.

\begin{prop}
\label{prop:so3-symplectic-action}
For each $v\in\mathfrak{so}(3)$, the corresponding infinitesimal generator $X_v$ has a Hamiltonian function $H_v:S^{2} \to \mathbb R$, where $\pi_{\Gamma_{S^2}}^* H_v = -f_{\bar\alpha}(\mathcal{X}_{v})$, and $\bar\alpha = f_{\bar\alpha} \, \frac{\partial}{\partial\theta}$ is an $\mathbb{SO}(3)$-invariant connection $1$-form.
\end{prop}

\begin{proof}
The $\mathbb{SO}(3)$-action on the principal $S^1$ bundle $\Gamma_{S^2}$ commutes with the inherent $S^{1}$-action of the bundle, and hence $\mathcal{X}_{v}$ is $S^{1}$-invariant. 
As a result, $f_{\bar\alpha}(\mathcal{X}_{v})$ is also $S^{1}$-invariant, and there is a function $H_{v}$ on $M$ such that $\pi_{\Gamma_{S^2}}^* H_v = -f_{\bar\alpha}(\mathcal{X}_{v})$.
By Cartan's formula,
\[
\begin{aligned}
\pi_{\Gamma_{S^2}}^* d H_v 
= -d\big(f_{\bar{\alpha}}(\mathcal{X}_{v})\big)
= \iota_{\mathcal{X}_{v}}df_{\bar\alpha} - \mathcal{L}_{\mathcal{X}_{v}}f_{\bar\alpha}
= \iota_{\mathcal{X}_{v}} \pi_{\Gamma_{S^2}}^* \omega_{S^2} 
= \pi_{\Gamma_{S^2}}^* \iota_{X_v} \omega_{S^2},
\end{aligned}
\]
where we used that $X_v = (\pi_{\Gamma_{S^2}})_* \mathcal{X}_v$ and that $f_{\bar\alpha}$ is $\mathbb{SO}(3)$-invariant and thus $\mathcal{X}_v$ invariant.
Since $(\pi_{\Gamma_{S^2}})_*$ is surjective, we obtain $dH_v = \iota_{X_v} \omega_{S^2}$, that is $X_v$ is Hamiltonian, with Hamiltonian function $H_v$.
\end{proof}

\section{Maslov bundles of symplectic homogeneous \texorpdfstring{$G$}{G}-spaces}
\label{sec:maslov-bundles-homogeneous}

In this subsection we extend Proposition \ref{prop:so3-transitive-maslov-bundle} to the case where $M$ is a symplectic homogeneous $G$-space, namely, the group action $\Phi$ by $G$ is transitive on $M$.
In particular, we show Theorem \ref{thm:homogeneous-G-space-(intro)}, that is, we show that when the first real Chern class $c_\Gamma$ of $\Gamma_J$ is nontrivial, then $\Gamma_J$ is a homogeneous $G$-space.

For each element $v$ of the Lie algebra $\mathfrak{g}$ of $G$, denote by $\mathcal{X}_{v}$ the infinitesimal generator of $\Phi_{\Gamma}$ in the direction $v$. 
That is, $\mathcal{X}_{v}(w)=\frac{d}{dt}\big|_{t=0} (\exp(tv)\cdot w)$ with $t\mapsto \exp(tv)$ being the one parameter subgroup of $G$ generated by $v$.

\begin{prop}
\label{prop:Phi-not-transitive-integrable-distribution}
If $\Phi_\Gamma$ is not transitive on $\Gamma_J$, then the infinitesimal generators $\mathcal{X}_{v}$, $v \in \mathfrak{g}$, span an integrable connection $\mathcal{D}$ on the bundle $\Gamma_{J}$.
\end{prop}

\begin{proof}
As an orbit of a compact Lie group action, $G\cdot w$ is an embedded
submanifold in $\Gamma_{J}$, and thus its dimension is no larger
than $\dim \Gamma_J = 2n+1$. 

If $\dim (G\cdot w) = 2n+1$, then $G\cdot w$ is an open set
in $\Gamma_{J}$ and $(G\cdot w) \cap\pi_{\Gamma_{J}}^{-1}(p)$
is open in $\pi_{\Gamma_{J}}^{-1}(p)$ and hence it is the whole $\pi_{\Gamma_{J}}^{-1}(p)$ since it is open and compact for any $p\in M$. 
This implies that the action $\Phi_\Gamma$ is transitive on $\Gamma_J$ thus contradicting the assumption of the proposition.

Now suppose that $\dim (G\cdot w) < 2n+1$. 
Since $G$ acts transitively on $M$ and $\pi_{\Gamma} \circ \Phi_{\#} = \Phi\circ\pi_{\Gamma}$, it holds
\[
(\pi_{\Gamma})_*\big(T_{w}(G\cdot w)\big)=T_{\pi_{\Gamma}(w)}M,
\]
and hence $\dim (G\cdot w) \ge 2n$. 
Together this yields $\dim (G\cdot w) = 2n$.
Due to the commutativity of $\Phi_{\Gamma}$ and the $S^{1}$ action on $\Gamma_{J}$, $(G\cdot w)\cdot z = G\cdot(w\cdot z)$ with $z\in S^1$.
This means $(G\cdot w)\cdot z$ is exactly the $G$-orbit through $w\cdot z$, and the $S^{1}$ action maps $G$-orbits to $G$-orbits. 
Since
\[
T_{w}(G\cdot w) = \operatorname{span} \{\mathcal{X}_{v}\big|v\in\mathfrak{g}\}=\mathcal{D}_{w},
\]
$\mathcal{D}$ is a $2n$ dimensional distribution invariant under the $S^{1}$ action, and $(\pi_{\Gamma})_*(\mathcal{D}_{w}) = T_{\pi_{\Gamma}(w)}M$.
Notice that the orbits $G\cdot w$, $w \in \Gamma_J$, are the maximal connected integral manifolds of $\mathcal{D}$.
\end{proof}

Using Proposition \ref{prop:Phi-not-transitive-integrable-distribution}, we can now prove Theorem~\ref{thm:homogeneous-G-space-(intro)} which we recall here.

\begin{thmiv}
Let $G$ be a compact Lie group acting transitively and symplectically on $M$, that is, $M$ is homogeneous $G$-space. 
If the Chern class $c_\Gamma$ of the Maslov $S^1$ bundle $\Gamma_J$ is non-zero, then $\Gamma_J$ is also a homogeneous $G$-space.
\end{thmiv}

\begin{proof}
If $\Phi_\Gamma$ does not act transitively on $\Gamma_J$, then according to Proposition \ref{prop:Phi-not-transitive-integrable-distribution}, there exists an integrable connection $\mathcal{D}$.
Then, Lemma~\ref{lem:characteristic-class-integrable-connection} implies that the Chern class of $\Gamma_J$ vanishes, contradicting the theorem's assumption.
\end{proof}

\section{Symplectic actions on monotone symplectic manifolds}
\label{sec:monotone-symplectic}

Consider a symplectic manifold $(M,\omega)$ and denote by $c_{\Gamma_J} \in H^2_{dR}(M)$ and $c_{\Gamma_J^2} \in H^2_{dR}(M)$ the real Chern classes of the Maslov $S^1$ bundles $\Gamma_J$ and $\Gamma_J^2$, respectively. Then, $c_{\Gamma_J^2} = 2 c_{\Gamma_J}$.

\begin{defn}
\label{def:monotone-symplectic}
A symplectic manifold $(M,\omega)$ is called \emph{monotone} if $[\omega] = r \, c_{\Gamma_J}$, for $r \in \mathbb{R}$.
\end{defn}

\begin{rem}
Definition~\ref{def:monotone-symplectic} differs from the definition of monotone symplectic manifold in \cite{Lupton1995} which requires that $r > 0$. 
\end{rem}

\begin{rem}
The monotone condition $[\omega] = r\, c_{\Gamma_J}$ implies $[\omega] = r'\,c_{\Gamma_J^2}$, with $r = 2r'$.
The arguments in this section apply to both $\Gamma_J$ or $\Gamma_J^2$ and below we write $\Gamma$ to refer to either one of the Maslov $S^1$ bundles and we write $c_\Gamma$ for the corresponding Chern class.
\end{rem}

Recall that a symplectic $G$-action on $M$ is Hamiltonian if and only if there is a smooth map
\[
\mu : M \to \mathfrak{g}^{*} : p \mapsto F_p,
\]
called \emph{momentum map}, such that for each $h\in G$ and $v\in\mathfrak{g}$, $\mu$ is $G$-coadjoint equivariant, i.e.,
\[
\mu_{\Phi^h(p)} = Ad_{h^{-1}}^{*} \mu_p,
\]
and for each $v \in \mathfrak{g}$, the mapping 
\[
H_v : M \to \mathbb{R} : p \mapsto H_v(p) = \mu_p(v)
\]
defines a Hamiltonian for the corresponding infinitesimal generator $X_v$, i.e., $dH_v = \iota_{X_v}\omega$.

\begin{defn}
An $S^{1}$-invariant $1$-form $\eta$ on the bundle $\pi_{\Gamma} : \Gamma \to M$ is called a \emph{symplectic potential} if $d\eta = -\pi_\Gamma^{*}\omega$.
\end{defn}

The following lemma asserts that monotone symplectic manifolds always have a symplectic potential.

\begin{lem}\label{lem:existence-symplectic-potential}
Let $(M,\omega)$ be a symplectic manifold, denote by $c_{\Gamma}\in H_{dR}^{2}(M)$ the first real Chern class of the Maslov $S^1$ bundle $\pi_\Gamma : \Gamma \to M$, and assume that $[\omega]=r \, c_{\Gamma}$ for some real number $r$.
Then for each connection $1$-form $\alpha=f_{\alpha}\,\frac{\partial}{\partial\theta}$ on $\Gamma$ there exists an $1$-form $\tau_\alpha$ on $M$ such that 
\begin{equation}
\label{eq:symplectic-potential-monotone}
\eta_\alpha = -r f_\alpha + \pi_\Gamma^*\tau_\alpha,
\end{equation}
is a symplectic potential.
\end{lem}

\begin{proof}
Given the connection $1$-form $\alpha$, denote by $\Omega_\alpha$ the corresponding curvature form, satisfying $df_\alpha = \pi_{\Gamma}^{*} \Omega_\alpha$.
Since $[\omega] = r \, c_\Gamma$ we have $[\omega - r \Omega_{\alpha}]= 0$, that is, there is an $1$-form $\tau_\alpha$ on $M$ such that $\omega = r \Omega_{\alpha}-d\tau_\alpha$. 
Let $\eta_\alpha = - r f_{\alpha} + \pi_\Gamma^*\tau_\alpha$. 
It is straightforward to check that $\eta_\alpha$ is invariant with respect to the inherent $S^1$ action on $\Gamma$.
Moreover,
\[ d\eta_\alpha
= -r \, df_\alpha + \pi_{\Gamma}^* d\tau_\alpha
= \pi_{\Gamma}^*(-r \, \Omega_\alpha + d\tau_\alpha) 
= -\pi_{\Gamma}^{*}\omega. \] 
Hence, $\eta_\alpha$ is a symplectic potential on the bundle $\pi_\Gamma: \Gamma \to M$.
\end{proof}

Averaging a symplectic potential $\eta_\alpha$ using the lifted $G$ action $\Phi_{\Gamma}$ gives a $G$-invariant symplectic potential $\bar\eta_\alpha$.
In particular, for $\eta_\alpha= -r f_\alpha + \pi_\Gamma^* \tau_\alpha$ in Eq.~\eqref{eq:symplectic-potential-monotone} we find
\[ 
\bar\eta_\alpha
= \int_G (\Phi_\Gamma^h)^* \eta_\alpha d\mu(h)
= -r \int_G (\Phi_\Gamma^h)^* f_\alpha d\mu(h)
+ \int_G (\Phi_\Gamma^h)^* \pi_\Gamma^* \tau_\alpha d\mu(h)
= -r \, \bar{f}_{\alpha} + \pi_\Gamma^* \bar\tau_\alpha,
\]
where $\bar\alpha = \bar f_\alpha \frac{\partial}{\partial\theta}$ is the result of averaging $\alpha$ over $G$.
Moreover, since the action $\Phi_\Gamma$ covers $\Phi$, that is, $\pi_\Gamma \circ \Phi_\Gamma^h = \Phi^h \circ \pi_\Gamma$, and $\Phi$ is a symplectic action, we have
\[ d \bar\eta_\alpha 
= \int_G (\Phi_\Gamma^h)^* d \eta_\alpha \, d\mu(h)
= -\int_G (\Phi_\Gamma^h)^* \pi_\Gamma^* \omega \, d\mu(h)
= -\int_G \pi_\Gamma^* (\Phi^h)^* \omega \, d\mu(h)
= -\pi_\Gamma^* \int_G \omega \, d\mu(h) = -\pi_\Gamma^* \omega.
\]

\begin{defn}
\label{def:associated-symplectic-potential}
Let $(M,\omega)$ be a monotone symplectic manifold with $[\omega]=r \, c_{\Gamma}$.
Given a connection $1$-form $\alpha=f_{\alpha} \frac{\partial}{\partial\theta}$ on $\Gamma$, the \emph{associated symplectic potential} is defined by $\eta_\alpha = -r f_\alpha + \pi_\Gamma^*\tau_\alpha$, where $d\tau_\alpha = r\Omega_\alpha - \omega$. If $\Phi$ is a symplectic action of a compact Lie group $G$ on $M$ and $\Phi_\Gamma$ is the lifted action on $\Gamma$, the \emph{associated $G$-invariant symplectic potential} $\bar\eta_\alpha$ is defined by averaging $\eta_\alpha$ over the action $\Phi_\Gamma$ of $G$.
\end{defn}

We have the following lemma.

\begin{lem}
\label{lem:coadjoint-equivariant map from a symplectic potential}
Let $\bar{\eta}$ be a $G$-invariant symplectic potential.
For each $p \in M$, let $\mu_p$ be the element in $\mathfrak{g}^{*}$ with $\mu_{p}(v)=\bar{\eta}(\mathcal{X}_{v})(w)$ for $v\in\mathfrak{g}$ and any $w \in (\pi_\Gamma)^{-1}(p)$, where $\mathcal{X}_{v}$ is the infinitesimal generator of $\Phi_{\Gamma}$ corresponding to $v$. 
Then, the map
\[
\mu: M \to \mathfrak{g}^* : p \mapsto \mu_{p}
\]
is a momentum map for $G$, and for each $v \in \mathfrak{g}$, the corresponding Hamiltonian function $H_v : M \to \mathbb{R}$ satisfies $\pi_\Gamma^* H_v = \bar{\eta}(\mathcal{X}_{v})$.
\end{lem}

\begin{proof}
Since $\bar{\eta}$ and $\mathcal{X}_{v}$ are both invariant with respect to the inherent $S^1$ action on $\Gamma$, the function $\bar{\eta}(\mathcal{X}_{v})$ is constant along each $S^1$ fiber of $\Gamma$, that is, $\bar{\eta}(\mathcal{X}_{v})(w) = \bar{\eta}(\mathcal{X}_{v})(w')$ for all $w, w' \in (\pi_\Gamma)^{-1}(p)$.
Therefore, $\mu$ is well defined.

For the smoothness of $\mu$, it suffices to show that the map $\bar{\mu} : \Gamma\times\mathfrak{g} \to \mathbb{R}$ defined by $\bar{\mu}(p,v) = \mu_{p}(v)$ is smooth. 
Observe that $\bar{\mu}$ factors as
\[
\bar{\mu}:\Gamma\times\mathfrak{g}\stackrel{\mathcal{X}}{\longrightarrow}T\Gamma\stackrel{\bar{\eta}}{\longrightarrow}\mathbb{R}
\]
with $\mathcal{X}$ being the map sending $(p,v)$ to $\mathcal{X}_{v}(p)$. 
The factorization below shows smoothness of $\mathcal{X}$:
\[
\mathcal{X}:\Gamma\times\mathfrak{g}\stackrel{\sigma}{\longrightarrow}T\Gamma\times TG\cong T(\Gamma\times G)\stackrel{D\Phi_{\Gamma}}{\longrightarrow}T\Gamma.
\]
Here, $\sigma$ is the map sending $(p,v)\in\Gamma\times\mathfrak{g}$ to $(\bar{\boldsymbol{0}}_{p},v)\in T_{p}\Gamma\times T_{1}G$ with
$\bar{\boldsymbol{0}}$ being the zero section from $\Gamma$ to $T\Gamma$,
and $D\Phi_{\Gamma}$ is the tangent map of $\Phi_{\Gamma}$.
Hence, $\mathcal{X}$ is smooth and thus $\bar{\mu}$ is also smooth.

To show that $\mu$ is $G$-coadjoint equivariant notice that since the action $\Phi_\Gamma$ of $G$ on $\Gamma$ covers the action $\Phi$ on $M$ we have
\[
\begin{aligned}
\mu_{\Phi^h(p)}(v) = \bar{\eta}(\mathcal{X}_{v})(\Phi_\Gamma^h(w))
= (\Phi_\Gamma^h)^* (\bar{\eta}(\mathcal{X}_{v})) (w).
\end{aligned}
\]
Since $\bar\eta$ is $\Phi_\Gamma$ invariant, we find
\[
(\Phi_\Gamma^h)^* (\bar{\eta}(\mathcal{X}_{v}))
= (\Phi_\Gamma^h)^* \bar{\eta} ((\Phi_\Gamma^h)^* (\mathcal{X}_v))
= \bar\eta (\mathcal{X}_{Ad_{h^{-1}}v}),
\]
where we used that for $h\in G$ and $v \in \mathfrak{g}$ we have $\mathcal{X}_{Ad_{h^{-1}}v} = (\Phi_\Gamma^h)^* \mathcal{X}_v$, see \cite{MarsdenRatiu1999}.
Therefore,
\[ 
\mu_{\Phi^h(p)}(v) = \bar\eta (\mathcal{X}_{Ad_{h^{-1}}v})(w)
= \mu_p(Ad_{h^{-1}}v),
\]
implying $\mu_{\Phi^h(p)}(v) = \mu_p \circ Ad_{h^{-1}} = Ad_{h^{-1}}^* \mu_p$.

For each $v \in \mathfrak{g}$ consider the function $H_v : M \to \mathbb R$ defined by $H_v(p) = \mu_p(v) = \bar\eta(\mathcal{X}_v)(w)$.
This implies $H_v \circ \pi_\Gamma (w) = \bar\eta(\mathcal{X}_v)(w)$ and thus $(\pi_\Gamma)^* H_v = \bar\eta(\mathcal{X}_v)$.
To show that $H_v$ is a Hamiltonian function for the infinitesimal generator $X_v$ of the action $\Phi$ on $M$ we notice that $dH_v = \iota_{X_v}\omega$ if and only if $d ((\pi_\Gamma)^* H_v) = (\pi_\Gamma)^* (\iota_{X_v}\omega)$ due to the surjectivity of $(\pi_\Gamma)_*$.
We check
\[ 
d ((\pi_\Gamma)^* H_v) 
= d(\bar\eta(\mathcal X_v)) 
= \mathcal{L}_{\mathcal X_v} \bar\eta - \iota_{\mathcal X_v} d\bar\eta
= \iota_{\mathcal X_v} ((\pi_\Gamma)^* \omega)
= (\pi_\Gamma)^* (\iota_{X_v} \omega),
\]
where we used that $X_v = (\pi_\Gamma)_* \mathcal X_v$. 
This concludes the proof.
\end{proof}

We can now prove Theorem \ref{thm:hamiltonian-G-action} which we restate here using the terminology introduced in this section.

\begin{thmv}\cite{Ono1988,Lupton1995}
If $(M,\omega)$ is a monotone symplectic manifold then any symplectic action $\Phi$ on $M$ by a compact Lie group $G$ is Hamiltonian. 
For any $v\in\mathfrak{g}$, there exists a Hamiltonian $H_v$ for the symplectic vector field $X_v$ on $M$, satisfying $\pi_{\Gamma}^*(H_v) = \bar{\eta}(\mathcal{X}_{v})$ for a $G$-invariant symplectic potential $\bar\eta$. 
\end{thmv}

\begin{proof}

Given a connection $1$-form $\alpha$ on $\Gamma$, Lemma~\ref{lem:existence-symplectic-potential} and the subsequent discussion imply the existence of an associated $G$-invariant symplectic potential $\bar{\eta}_\alpha$ satisfying Eq.~\eqref{eq:symplectic-potential-monotone}.
Then, by Lemma \ref{lem:coadjoint-equivariant map from a symplectic potential}, the $G$-action $\Phi$ is Hamiltonian and for $v \in \mathfrak{g}$ the corresponding Hamiltonian function $H_v$ satisfies $\pi_\Gamma^* H_v = \bar{\eta}_\alpha(\mathcal{X}_v)$.
\end{proof}

\section{Maslov data and local Maslov index for symplectic \texorpdfstring{$S^1$}{S1} actions}
\label{sec:maslov-data-s1}

In this section we consider the case where $G=S^{1}$ acts symplectically on the manifold $(M,\omega)$.
We denote by $\frac{\partial}{\partial\phi}$ the generator of the Lie algebra $\mathfrak{g} \cong \mathbb R \cdot \frac{\partial}{\partial\phi}$ of $G$, by $X_\phi$ the infinitesimal generator of the $S^1$ action $\Phi$ of $G$ on $M$, and by $\mathcal{X}_\phi$ the infinitesimal generator of the lifted $S^1$ action $\Phi_{\Gamma^2}$ on $\Gamma_J^2$.
We refer to the $G$ action on $\Gamma_J^2$ as the $S^1_\phi$ action, and we refer to the inherent $S^1$ action on $\Gamma_J^2$ as the $S^1_\theta$ action, to better distinguish between these two $S^1$ actions on $\Gamma_J^2$.
We denote by $S^1_\phi \cdot p$ the orbit of the $S^1$ action $\Phi$ on $M$ going through $p \in M$, and by $S^1_\phi \cdot w$ the orbit of the $S^1$ action $\Phi_\Gamma$ on $\Gamma_J^2$ going through $w \in \Gamma_J^2$. 
Finally, we denote by $S^1_\theta \cdot w$ the orbit of the inherent $S^1$ action on $\Gamma_J^2$ going through $w \in \Gamma_J^2$, that is, $S^1_\theta \cdot w = \pi_{\Gamma_J^2}^{-1}(\pi_{\Gamma_J^2}(w))$. 

Recall from Definition~\ref{def:maslov-data} that, given a connection $1$-form $\alpha = f_\alpha \, \frac{\partial}{\partial\theta}$ on $\Gamma_J^2$, the Maslov data of a smooth loop $\gamma : S^1 \to \Gamma_J^2$ with respect to $\alpha$ is 
\[ \mathfrak{md}_\alpha(\gamma) = \int_\gamma f_\alpha. \]

\begin{defn}
Consider a symplectic $S^1$ action $\Phi$ on $(M,\omega)$. 
Then the \emph{Maslov data of $\Phi$ at $p \in M$ with respect to the connection $1$-form $\alpha$} is defined by
\begin{equation}
\label{eq:def-maslov-data-s1}
Q_{\alpha}(p)
= \mathfrak{md}_{\alpha}(\gamma_{w})
= \int_{\gamma_{w}} f_\alpha,
\end{equation}
where $\gamma_w = S_\phi^1 \cdot w$ for any $w\in\pi_{\Gamma_J^2}^{-1}(p)$. 
\end{defn}

The integral in Eq.~\eqref{eq:def-maslov-data-s1} is independent of the choice of $w$ on $\pi_{\Gamma_J^{2}}^{-1}(p)$ and hence $Q_{\alpha}$ is a well defined function on $M$. 
To see this, note that for $w,w'\in\pi_{\Gamma_J^{2}}^{-1}(p)$, there exists an element $u \in S^1_\theta$ such that $w' = w \cdot u$ and thus $\gamma_{w'} = \gamma_{w} \cdot u$, since $S^1_\phi$ and $S^1_\theta$ commute.
Then, $\int_{\gamma_{w}} f_\alpha = \int_{\gamma_{w'}} f_\alpha$ follows from the fact that $\alpha$ is a connection $1$-form and thus $f_\alpha$ is invariant under $S^1_\theta$.

We now consider the Maslov data $Q_\alpha(p)$ at a fixed point $p$ of the $S^1$ action on $M$. 

\begin{defn}
\label{def:local-maslov-index-s1}
Given a fixed point $p$ of a symplectic $S^1$ action on $(M,\omega)$, the \emph{local Maslov index of the $S^1$ action at $p$} is defined as $k_p = Q_\alpha(p)$, where $\alpha$ is any connection $1$-form for the bundle $\Gamma_J^2$.
\end{defn}

Even though the definition of the local Maslov index $k_p$ makes use of a connection $1$-form $\alpha$, it turns out that its value is independent of the choice of $\alpha$. In particular, we have the following result.

\begin{prop}
\label{prop:local-maslov-index-value}
The local Maslov index $k_p = Q_\alpha(p)$ at a fixed point $p$ of a symplectic $S^1$ action on $(M,\omega)$ is an integer measuring how many times the orbit $S_\phi^1 \cdot w$ winds around the fiber $\pi_{\Gamma_J^2}^{-1}(p)$ and it does not depend on the choice of $\alpha$. 
\end{prop}

\begin{proof}
Suppose that $p$ is a fixed point of the action $\Phi$.
Then
\begin{equation} 
\pi_{\Gamma_J^{2}}^{-1}(p) = S^1_\theta \cdot w' = S^1_\phi \cdot w, 
\label{eq:fibre-orbit-fixed-point}
\end{equation}
for any $w, w' \in \pi_{\Gamma_J^{2}}^{-1}(p)$.
Eq.~\eqref{eq:fibre-orbit-fixed-point} also implies that for $w \in \pi_{\Gamma_J^{2}}^{-1}(p)$ we have $\mathcal{X}_\phi(w) = a(w) \frac{\partial}{\partial\theta}$.
Since the lifted action $S^1_\phi$ commutes with $S^1_\theta$, $\mathcal{X}_\phi$ is constant on $\pi_{\Gamma_J^{2}}^{-1}(p)$,
that is, for all $w\in\pi_{\Gamma_J^{2}}^{-1}(p)$ we have $\mathcal{X}_\phi(w) = a_p \frac{\partial}{\partial\theta}$, for some constant $a_p$, and thus $f_\alpha(\mathcal{X}_\phi) = a_p$ for any connection $1$-form $\alpha$ on $\Gamma_J^2$.
Hence,
\begin{equation*}
k_p 
= Q_{\alpha}(p)
= \int_{\gamma_w} f_\alpha
= \int_{0}^{1} f_\alpha(\mathcal{X}_\phi)_{\gamma_{w}(t)} dt
= a_p,
\end{equation*}
does not depend on $\alpha$.
The orbit of $\mathcal{X}_\phi = k_p \frac{\partial}{\partial\theta}$ is a closed orbit of $S^1_\phi$ with time $1$ recurrence. 
Therefore, $|k_p|$ counts the number of times the $S^1_\phi$ orbit winds around $\pi_{\Gamma_J^{2}}^{-1}(p)$.
\qedhere 
\end{proof}

The local Maslov index of a fixed point $p$ can be expressed in terms of the weights of the linearized $S^1$ action at $p$. In particular, if $p$ is a fixed point of the symplectic $S^1$ action $\Phi$, then a Darboux chart
\[
\Bigl(U, \sum_j dx_{j}\wedge dy_{j} \Bigr) \xrightarrow{\varphi} (M,\omega)
\]
with $\varphi(0)=p$ can be chosen in some $S^1$ invariant neighbourhood of $p$,
such that in this chart the restricted $S^1$ action is linearized as 
\[ e^{2 \pi i t}, (z_1, \dots, z_n)  \mapsto (z_1 e^{2m_{1}\pi i t},\dots, z_n e^{2m_{n}\pi i t}), \quad z_j = y_j + i x_j. \]
We denote this linearized $S^{1}$ action on $U$ by $L_{\Phi}$
and we call $(m_{1},...,m_{n}) \in \mathbb{Z}^{n}$ the resonance type of
the fixed point $p$. 

\begin{prop}
\label{prop:Maslov data at the fixed points}
At a fixed point $p$ with resonance type $(m_{1},...,m_{n})$, the local Maslov index is
\begin{equation}
  \label{eq:maslov-data-resonant-type}
  k_{p} = 2 \sum_{i=1}^{n}m_{i}.
\end{equation}
\end{prop}

\begin{proof}
Denote by $\Gamma_{U}^{2}$ the Maslov $S^{1}$ bundle for $(U, \sum_j dx_j\wedge dy_j)$ and note that it is isomorphic to $U\times S^{1}$.
Since $\varphi$ is a symplectomorphism, it induces a bundle isomorphism $\varphi_*$ from $\Gamma_{U}^{2}$ to $\Gamma_{J}^{2}\big|_{\tilde{U}}$ with $\tilde{U}=\varphi(U)$. 
The lifted $S^{1}$ action $L_{\Phi}$ on $\Gamma_{U}^{2}$ and the action $\Phi_{\Gamma^{2}}$ on $\Gamma_{J}^{2}\big|_{\tilde{U}}$ are topologically conjugate via $\varphi_{*}$. 
Then it is straightforward to check that the local Maslov index of $L_{\Phi}$ at $0\in U$ is the same as that of $\Phi_{\Gamma^{2}}$ at $p$ and they both equal $k_{p}$.
For $(z_1,\dots,z_n,w) \in U \times S^1$, the lifted $S^1$ action $L_\Phi$ is given by
\[
e^{2\pi i t}, (z_1,\dots,z_n,w) \mapsto 
( e^{2\pi i m_1 t} z_1, \dots, e^{2\pi i m_n t} z_n, e^{4\pi i t \sum_{i=1}^n m_i} w).
\]
Therefore, the number of times that the $L_\Phi$ orbit winds around the fiber of $\Gamma_U^2$ at $0$, equals $2\sum_{i=1}^{n}m_{i}$, and it follows from Proposition~\ref{prop:local-maslov-index-value} that this is the local Maslov index $k_p$.
\end{proof}

In the remaining part of this section we discuss Maslov data and the local Maslov index in the specific cases of monotone symplectic manifolds, Maslov bundles with zero Chern class, trivial Maslov bundles, and cotangent bundles.

\subsection{Monotone symplectic manifolds}

We first turn our attention to monotone symplectic manifolds with $[\omega] = r c_{\Gamma^2}$, where $c_{\Gamma^2} \in H^2_{dR}(M)$ is the Chern class of the bundle $\Gamma_J^2$. 
Recall from Theorem~\ref{thm:hamiltonian-G-action} that in this case the vector field $X_\phi$ is Hamiltonian. 
We denote the corresponding Hamiltonian function by $H_\phi$.
Consider a $S^1_\phi$ invariant connection $1$-form $\alpha$ on $\Gamma_J^2$ and the associated $S^1_\phi$ invariant symplectic potential $\bar\eta_\alpha = -r f_{\bar\alpha} + \pi_{\Gamma_J^2}^* \bar\tau_\alpha$, see Definition~\ref{def:associated-symplectic-potential}.
Then we have $\pi_{\Gamma_J^2}^* H_\phi = \bar\eta_\alpha(\mathcal{X}_\phi)$.

\begin{lem}
\label{prop:H-phi}
Suppose that $S^1$ acts symplectically on the monotone symplectic manifold $(M,\omega)$, with $[\omega] = r \, c_{\Gamma^2}$.
Then the Hamiltonian $H_\phi$ satisfies
\begin{equation}\label{eq:H-phi}
H_\phi = -r Q_{\bar\alpha} + \bar\tau_\alpha(X_\phi). 
\end{equation}
\end{lem}

\begin{proof}
We have 
\[ \pi_{\Gamma_J^2}^* H_\phi 
= \bar\eta_\alpha(\mathcal X_\phi) 
= -r f_{\bar\alpha}(\mathcal X_\phi)
  + \pi_{\Gamma_J^2}^* \bar\tau_\alpha(X_\phi).
\]
Since $f_{\bar\alpha}(\mathcal X_\phi)$ is constant along $\gamma_w$ we obtain
\[ \pi_{\Gamma_J^2}^* Q_{\bar\alpha} 
= \int_{\gamma_w} f_{\bar\alpha}
= \int_0^1 f_{\bar\alpha}(\mathcal X_\phi)_{\gamma_w(t)} dt
= f_{\bar\alpha}(\mathcal X_\phi). \]
Therefore,
\[ \pi_{\Gamma_J^2}^* H_\phi 
= -r \pi_{\Gamma_J^2}^* Q_{\bar\alpha}
  + \pi_{\Gamma_J^2}^* \bar\tau_\alpha(X_\phi),
\]
and Eq.~\eqref{eq:H-phi} follows from the surjectivity of $\pi_{\Gamma_J^2}$.
\end{proof}

The next result follows directly from Proposition~\ref{prop:Maslov data at the fixed points} and Lemma~\ref{prop:H-phi}, using that if $p$ is a fixed point of the $S_\phi^1$ action then $Q_{\bar\alpha}(p) = k_p$ and $X_\phi(p) = 0$.

\begin{prop}\label{cor:H-phi-value-fixed-points}
Suppose that $S^1$ acts symplectically on the monotone symplectic manifold $(M,\omega)$, with $[\omega] = r \, c_{\Gamma^2}$. If $p$ is a fixed point of the $S^1$ action on $M$ with local Maslov index $k_p$, then $p$ is a critical point of $H_\phi$, and the critical value $H_\phi(p) = - r k_p \in r \cdot \mathbb Z$ depends only on the resonance type of $p$.
\end{prop}

\subsection{Maslov bundles with zero Chern class}

Proposition \ref{prop:local-maslov-index-value} shows that although the Maslov data $Q_{\alpha}(p)$ at arbitrary $p \in M$ depends on the choice of the connection $1$-form $\alpha$, its values at the fixed points of the $S^1$ action on $M$ do not have such dependence.
This implies that if the $S^1$ action has fixed points with different local Maslov indices, then there is no connection $1$-form $\alpha$ such that $Q_{\alpha}$ is constant on $M$. In this section we show that this may occur only when the Maslov bundle has a non-zero Chern class. In particular, we show the following result concerning Maslov bundles for which the Chern class $c_{\Gamma^2}$ vanishes.

\begin{prop}
\label{prop:vanishing-chern-maslov-data}
Suppose that $(M,\omega)$ is a connected symplectic manifold with $c_{\Gamma^2}=0$, and $\Phi$ is a symplectic $S^1$ action on $M$. 
Then there is a connection $1$-form $\alpha_0$ such that $Q_{\alpha_0}$ is constant on $M$. 
\end{prop}

\begin{proof}
Since $c_{\Gamma^2}=0$, if $\Omega_\alpha$ is the curvature form corresponding to an arbitrary connection $1$-form $\alpha = f_\alpha \frac{\partial}{\partial\theta}$, then there is an $1$-form $\tau$ on $M$ such that $\Omega_\alpha = d\tau$. 
Then, $\alpha_0 = f_{\alpha_0} \frac{\partial}{\partial\theta}$ with $f_{\alpha_0} = f_\alpha - \pi_{\Gamma_J^2}^* \tau$ is a connection $1$-form.
Moreover, we have
\[ df_{\alpha_0} = df_\alpha - \pi_{\Gamma_J^2}^* d\tau 
= \pi_{\Gamma_J^2}^* (\Omega_\alpha - d\tau) = 0. \]
For any two points $p_0, p_1 \in M$, consider corresponding points $w_0 \in (\pi_{\Gamma^2})^{-1}(p_0)$ and $w_1 \in (\pi_{\Gamma^2})^{-1}(p_1)$.
Let $\lambda:[0,1]\rightarrow\Gamma_{J}^{2}$
be a smooth path from $w_{0}$ to $w_{1}$, and
define the map $h : S^1\times[0,1] \to \Gamma_J^2$ by
\[
h(z,t)=\Phi_{\Gamma^{2}}^{z}\circ\lambda(t).
\]
Let $\tilde{f}_{\alpha_0}$ be the pullback of $f_{\alpha_0}$ to $S^{1}\times[0,1]$, i.e., $\tilde{f}_{\alpha_0}=h^* f_{\alpha_0}$.
Then $d\tilde{f}_{\alpha_0} = h^* df_{\alpha_0} = 0$. 
By Stokes' formula,
\begin{equation}
\int_{S^1\times\{0\}} \tilde{f}_{\alpha_0}
= \int_{S^1\times\{1\}} \tilde{f}_{\alpha_0}.
\label{eq:boundary integration}
\end{equation}
Note that for $i=0,1$, it holds $h_*(\frac{\partial}{\partial\theta}\big|_{(z,i)})
= \mathcal{X}_{\phi}\big(\Phi_{\Gamma^{2}}^{z}(w_{i})\big)$.
As a consequence, we have
\[
\int_{S^1\times\{i\}} \tilde{f}_{\alpha_0}
= \int_0^1 \tilde{f}_{\alpha_0} \Bigl(\frac{\partial}{\partial\theta} \Bigr) \Big|_{(e^{2\pi i t},i)} dt
= \int_0^1 f_{\alpha_0}(\mathcal{X}_{\phi}) \big|_{\gamma_{w_{i}}(t)} dt
= \int_{\gamma_{w_{i}}} f_{\alpha} = Q_{\alpha_0}(p_i),
\]
where $\gamma_{w_{i}}(t)=\Phi_{\Gamma^{2}}^{e^{2\pi i t}}(w_{i})$. Since $p_0$, $p_1$ are arbitrary points in $M$, we obtain that $Q_{\alpha_0}(p)$ is constant on $M$.
\end{proof}

The following result is a direct consequence of Proposition~\ref{prop:local-maslov-index-value} and Proposition~\ref{prop:vanishing-chern-maslov-data}.

\begin{cor}
If $(M,\omega)$ is a connected symplectic manifold with $c_{\Gamma^2}=0$, and $\Phi$ is a symplectic $S^1$ action on $M$, then all fixed points of $\Phi$ have the same local Maslov index.
\end{cor}

\subsection{Trivial Maslov bundles}

We consider now the case where $\Gamma_J^{2}$ is a trivial bundle and thus the Maslov index can be defined, see Definition~\ref{def:maslov-index(intro)}.
Notice that in this case we also have $c_{\Gamma^2} = 0$.
Let $\mathfrak{s}$ be a global section of $\Gamma_J^2$.
If $L \subset T_p M$ is a Lagrangian plane, then $\tilde{\gamma}_{L}(z) = \Phi_{*}^{z}(L)$ with $z \in S^1$ is a loop in $\Lambda_{J}$.
Hence, its Maslov index with respect to $\mathfrak{s}$, denoted by $\mathfrak{m}_{\mathfrak{s}}(\tilde{\gamma}_{L})$, is the degree of the map obtained by the composition
\begin{equation}
S^1
\xrightarrow{\tilde\gamma_L} \Lambda_{J}
\xrightarrow{det_{J}^{2}}\Gamma_{J}^{2}
\xrightarrow{tr_{\mathfrak{s}}}
M\times S^{1}\xrightarrow{pr_{S^{1}}}S^{1}.
\label{eq:Maslov index for circle actions (Lambda)}
\end{equation}
Equivalently, we can consider the degree $\mathfrak{m}_{\mathfrak{s}}(\gamma_w)$ of the map
\begin{equation}
S^1
\xrightarrow{\gamma_w}\Gamma_{J}^{2}
\xrightarrow{tr_{\mathfrak{s}}}M\times S^{1}
\xrightarrow{pr_{S^{1}}}S^{1}
\label{eq:Maslov index for circle actions (Gamma)}
\end{equation}
with $\gamma_{w}(z)=\Phi_{\Gamma^{2}}^{z}(w)$ for $z \in S^1$.
The connectedness of $\Gamma_J^2$, implies that $\mathfrak{m}_{\mathfrak{s}}(\gamma_w)$ is independent of $w \in \Gamma_J^2$ and thus all $S^1_\phi$ orbits $\gamma_w$ have the same Maslov index. 

\begin{defn}
Consider a symplectic manifold $M$ such that $\Gamma_J^2$ is trivial, and a symplectic $S^1$ action $\Phi$ on $M$. 
The \emph{Maslov index of the $S^1$ action $\Phi$} with respect to a section $\mathfrak{s}:M \to \Gamma_J^2$ is defined as the Maslov index $\mathfrak{m}_{\mathfrak{s}}(\gamma_w)$ for any $w \in \Gamma_J^2$ and is denoted by $\mathfrak{m}_{\mathfrak{s}}(\Phi)$.
\end{defn}

If $p \in M$ is a fixed point of the $S^1$ action, then the local Maslov index $k_p$ is the degree of the map in Eq.~\eqref{eq:Maslov index for circle actions (Gamma)}, that is, $k_p = \mathfrak{m}_{\mathfrak{s}}(\Phi)$.
Since $k_p$ does not depend on the choice of section $\mathfrak{s}$ we conclude that $\mathfrak{m}_{\mathfrak{s}}(\Phi)$ also does not depend on $\mathfrak{s}$.
Therefore, we obtain the following statement.

\begin{prop}
\label{prop:maslov-index-phi}
When the Maslov $S^1$ bundle $\Gamma_J^2$ of a connected manifold $M$ is trivial and the $S^1$ action $\Phi$ on $M$ has fixed points, then the Maslov index $\mathfrak{m}_{\mathfrak{s}}(\Phi)$ of $\Phi$ does not depend on the choice of section $\mathfrak{s} : M \to \Gamma_J^2$ and $\mathfrak{m}_{\mathfrak{s}}(\Phi) = k_p$, where $p \in M$ is any fixed point of $\Phi$.
\end{prop}

Note that if the $S^1$ action $\Phi$ does not have fixed points, the Maslov index of the $S^1$ action may depend on the choice of section $\mathfrak{s}$.
An example where this occurs is given by the action of $S^1$ on the cylinder $\mathbb{C}^*:=\mathbb{C}\setminus\{0\}$ with the standard symplectic form $dx \wedge dy$, given by $(e^{i\phi}, z ) \in S^1 \times \mathbb{C}^* \mapsto e^{i\phi}z \in \mathbb{C}^*$.
Then, $\Gamma_J^2 \simeq \mathbb C^* \times S^1$ is trivial.
The vectors $\frac{\partial}{\partial x}$ define a Lagrangian distribution with respect to which the Maslov index of the $S^1$ orbits is $2$, while the vectors $\frac{\partial}{\partial \phi}$ define a Lagrangian distribution with respect to which the Maslov index of the $S^1$ orbits is $0$.

Finally, suppose that $\mathcal{S} \subset M$ is a Lagrangian submanifold invariant under the symplectic $S^1$ action $\Phi$.
Then, given a point $p \in \mathcal{S}$, $\gamma(z) = \Phi^z(p)$ for $z \in S^1$ is a loop on $\mathcal{S}$ and its Maslov index $\mathfrak{m}_\mathfrak{s}(\gamma)$ is the degree of the map
\[
S^1
\xrightarrow{\sigma_{\mathcal{S}} \circ \gamma}
\Lambda_{J}
\stackrel{det_J^2}{\longrightarrow}
\Gamma_{J}^{2}\stackrel{tr_{\mathfrak{s}}}{\longrightarrow}
M\times S^{1}
\stackrel{pr_{S^1}}{\longrightarrow}
S^1,
\]
where $\sigma_{\mathcal{S}} \circ \gamma(z) = T_{\gamma(z)} \mathcal{S}$, see Definition~\ref{def:maslov-index-lagrangian}.
Since $T_{\gamma(z)} \mathcal{S} = T_{\Phi^{z}(p)} \mathcal{S} = \Phi_*^z(T_p\mathcal{S}) = \tilde{\gamma}_{T_p \mathcal{S}}$,
the last map is exactly the map in Eq.~\eqref{eq:Maslov index for circle actions (Lambda)} for $L = T_p \mathcal{S}$, showing that in this case the Maslov index of the $S^1$ action $\Phi$ coincides with the usual Maslov index $\mathfrak{m}_\mathfrak{s}(\gamma)$.

\subsection{Cotangent bundles}
\label{sec:cotangent-bundles}

According to \cite{Ono1988, Lupton1995, Cho2017}, when $(M,\omega)$ is compact and $c_{\Gamma} = 0$, there is no effective Hamiltonian circle action. 
This result, however, does not apply to the case which is of the most concern for physics, that is, when $M$ is a cotangent bundle. 
In this case $c_{\Gamma} = [\omega] = 0$, and the Maslov bundles are trivial. 

To understand the difference between compact and non-compact phase spaces we briefly review the approach used in \cite{Cho2017} to show that there is no effective Hamiltonian circle action in the compact case. 
The proof in \cite{Cho2017} relies on the construction of an embedded $2$-sphere by taking the $S^1$-action on an orbit of the gradient flow of a supposed Hamiltonian. 
This construction is possible when $M$ is compact since in this case any nontrivial orbit of the gradient flow goes forward and backward to different fixed points of the $S^1$ action. 
However, in the case where $M$ is non-compact, an orbit of the gradient flow may not have both a forward limit and a backward limit. 
Note that this is compatible with Remark~\ref{rem:unique critical value of S^1 Hamiltonian on cotangent bundle}: since there is at most one critical value of the Hamiltonian and the value of the Hamiltonian strictly increases along a nontrivial orbit of the gradient flow, no orbit can have both a forward limit and a backward limit.

We conclude this discussion on cotangent bundles with the following two remarks which can be obtained by averaging the Liouville $1$-form of a cotangent bundle, but they also appear as consequences of the discussion in this section about local Maslov indices. 

\begin{rem}\label{rem:unique critical value of S^1 Hamiltonian on cotangent bundle}
Suppose that $M$ is a connected cotangent bundle and the symplectic $S^1$ action on $M$ is Hamiltonian with Hamiltonian function $H_\phi$.
Since all fixed points of the $S^1$ action have the same local Maslov index, $H_\phi$ has at most $1$ critical value.
\end{rem}

\begin{rem}
Consider the case where $M$ is a cotangent bundle and $\dim M=4$.
The symplectic $S^{1}$ action can be linearized in a Darboux chart in a neighbourhood of a fixed point $p$ as $t\mapsto(e^{2\pi i m t},e^{-2\pi i nt})$.
When $k_p=0$, we have $m-n=0$, and hence all the fixed points have resonance types $(1,-1)$.
\end{rem}

\section{Applications to integrable Hamiltonian systems}
\label{sec:int-ham-sys}

An integrable Hamiltonian system with $2$ degrees of freedom is a triple $(M,\Phi,F)$ with a $4$ dimensional symplectic manifold $M$, a Hamiltonian $\mathbb{R}^{2}$ action $\Phi$, and an integral map $F : M \to \mathbb{R}^2$ such that $\Phi^* F = F$, the Poisson bracket $\{F_1, F_2\} = 0$ vanishes, and $dF_1 \wedge dF_2 \ne 0$ almost everywhere.
The regular domain $\mathcal{D}_\mathrm{reg}$ of the system is the part of the phase space $M$ that consists of compact regular orbits of $\Phi$.

In this section we discuss two applications of the concepts and the results on local Maslov indices of $S^1$ actions from Sec.~\ref{sec:maslov-data-s1} to integrable Hamiltonian systems.
First, in Sec.~\ref{sec:additional-s1}, we discuss necessary conditions under which, if an integrable Hamiltonian system possesses an $S^1$ action in $\Dreg$, this can be extended to a $\mathbb T^2$ action.
Then, in Sec.~\ref{sec:pinched-tori}, we show that the Maslov $S^1$ bundle over a pinched torus is trivial.

\subsection{Extending an \texorpdfstring{$S^1$}{S1} action to a \texorpdfstring{$\mathbb T^2$}{T2} action}
\label{sec:additional-s1}

Several physically important integrable Hamiltonian systems, such as the spherical pendulum, the isotropic planar harmonic oscillator, and coupled spin systems, have a Hamiltonian $S^1$ symmetry. 
That is, in these systems, there is a first integral generating an $S^1$ action. 
Integrable Hamiltonian systems with global action-angle coordinates have a Hamiltonian $\mathbb T^2$ action. 
From this point of view, it is interesting to understand, if a system with an $S^1$ action can have an additional, independent, $S^1$ action, giving rise to a $\mathbb T^2$ action for the system. 
We have the following result.

\begin{prop}
\label{prop:action-angle coordinates and local Maslov indices}
Let $(M, \Phi, F)$ be an integrable Hamiltonian system with $2$ degrees of freedom, where the integral map $F$ is proper and $\Gamma_J^2 \cong M\times S^1$. 
Suppose that the system $(M, \Phi, F)$ has an $S^1$ action $\varphi_\theta$ with nontrivial local Maslov index at its fixed points. 
Then the system possesses an additional, independent, $S^1$ action on $\mathcal{D}_\mathrm{reg}$.
\end{prop}

\begin{proof}
By definition, the regular domain $\mathcal{D}_\mathrm{reg}$ is foliated by compact orbits of the Hamiltonian $\mathbb R^2$ action $\Phi$.
These compact orbits are homeomorphic to $\mathbb T^2$. 
Let $\mathcal{O}$ be the orbit space of the $\mathbb R^2$ action in $\mathcal{D}_\mathrm{reg}$, and for each $o \in \mathcal{O}$ denote by $\mathcal{L}_o$ the corresponding period lattice which is isomorphic to $\mathbb Z^2$.
The period lattices of the orbits form the smooth period lattice bundle $\mathcal L \to \mathcal O $, where $\mathcal{L}$ is the disjoint union of the period lattices $\mathcal{L}_o$, $o \in \mathcal{O}$.
    
Fix a trivialization $\Gamma^2_J \cong M \times S^1$ by choosing a global section $\mathfrak{s}: M \to \Gamma_J^2$. 
For each $(o,l) \in \mathcal{L}$ consider the corresponding loop $\gamma_l$ on  $o$, and denote its Maslov index by $\mathfrak{m_s}(\gamma_l)$.
This defines the bundle homomorphism
\[ \rho: \mathcal{L} \to \mathbb{Z} 
   : (o,l) \mapsto \rho(o,l) = \mathfrak{m_s}(\gamma_l). \]

For any path in $\mathcal O$ from $o$ to $o'$, the parallel transport of points in $\mathcal L_o$ along the path is well defined and gives an isomorphism $\mu : \mathcal L_o \to \mathcal L_{o'}$. 
Since $\rho$ takes values in $\mathbb Z$, its value remains constant under parallel transport and thus $\rho|_{\mathcal L_{o'}} \circ \mu = \rho|_{\mathcal L_o}$.

For each $o \in \mathcal O$ choose $\gamma$ to be an orbit of the $S^1$ action $\varphi_\theta$ on $o$, and denote by $l_1(o) \in \mathcal L_o$ the corresponding element of the period lattice. 
Since the $S^1$ action is defined on $\mathcal D_\mathrm{reg}$, the mapping  $o \mapsto l_1(o)$ is a global section of the period lattice bundle $\mathcal L \to \mathcal O$. 
In particular, the parallel transport of $l_1(o)$ along any path from $o$ to $o'$ equals $\mu(l_1(o)) = l_1(o')$.
The Maslov index $\rho(o,l_1(o)) = \mathfrak{m}_\mathfrak{s}(\gamma)$ is the Maslov index of the $S^1$ action, i.e., $\rho(o,l_1(o)) =\mathfrak{m}_\mathfrak{s}(\varphi_\theta)$. 
Since $\varphi_\theta$ has fixed points, $\mathfrak{m}_\mathfrak{s}(\varphi_\theta)$ equals the local Maslov index at these points (Proposition \ref{prop:maslov-index-phi}) and is thus, by assumption, not zero. 
Therefore, the homomorphism $\rho|_{\mathcal L_o} : \mathcal L_o \to \mathbb Z$ is nontrivial and $\ker \rho|_{\mathcal L_o}$ is one dimensional. 

Denote by $l_* \in \mathcal L_{o_*}$ a generator of $\ker \rho|_{\mathcal L_{o_*}}$ over $\mathbb Z$ for a fixed $o_* \in \mathcal O$.
If $\mu$ is the parallel transport from the fixed $o_*$ to an arbitrary $o \in \mathcal O$, then $\rho(\mu(l_*)) = \rho(l_*) = 0$, that is, $\mu(l_*) \in \ker \rho|_{\mathcal L_{o}}$.
Moreover, $\mu$ induces an isomorphism from $\ker \rho|_{\mathcal L_{o_*}}$ to $\ker \rho|_{\mathcal L_{o}}$, and thus $\mu(l_*)$ is a generator of $\ker \rho|_{\mathcal L_{o}}$.
Consider the map $l_2: \mathcal O \to \mathcal L$ given by $l_2(o) = \mu(l_*)$ for $o \in \mathcal O$.
For any loop based at $o_*$ we have $\mu(l_*) = l_*$ since parallel transport along any loop based at $o$ preserves the orientation of $\mathcal L_{o_*}$.
This implies that the map $l_2$ is well defined since it is independent of the path from $o_*$ to $o$.
Therefore, $l_2 : \mathcal O \to \mathcal L$ defines a second, independent, global section of $\mathcal L \to \mathcal O$, and corresponds to an additional $S^1$ action on $\mathcal{D}_\mathrm{reg}$. 
\end{proof}

\subsection{Maslov \texorpdfstring{$S^1$}{S1} bundles over Lagrangian pinched tori}
\label{sec:Maslov S^1 Bundles over Lagrangian Pinched Tori}
\label{sec:pinched-tori}

Pinched tori appear as singular fibers in integrable Hamiltonian systems of $2$ degrees of freedom and they contain one or more focus-focus equilibria. 
In this discussion we consider only singly pinched tori; 
they consist of a fixed point $p_0$ of the Hamiltonian $\mathbb R^2$ action $\Phi$, and a regular orbit $o$ of $\Phi$ which is a Lagrangian submanifold diffeomorphic to $S^1 \times \mathbb{R}$. 
A pinched torus is homeomorphic to the quotient space $\big(S^{1}\times[-1,1]\big)\big/\big(S^{1}\times\{\pm1\}\big)$
which identifies the ends of the cylinder $S^{1}\times[-1,1]$ with a single point, or, equivalently, the quotient space $S^{2}\big/\{(0,0,\pm1)\}$ which glues the point $(0,0,1)$ with $(0,0,-1)$. 
Pinched tori in integrable Hamiltonian systems are associated to the appearance of Hamiltonian monodromy \cite{Duistermaat1980, Cushman1997, Zung1997}, but they also appear in a non-Hamiltonian context \cite{Cushman2001}. 

Here we focus on integrable Hamiltonian systems and consider pinched tori that are Lagrangian submanifolds of the symplectic manifold $(M,\omega)$. 
A Lagrangian pinched torus in $M$ can be represented by a continuous map
\[
\mathfrak{p}:S^{1}\times[-1,1]\rightarrow M,
\]
where the restriction $\mathfrak{p}|_{S^1 \times (-1,1)}$ is a Lagrangian embedding, and $\mathfrak{p}$ factors as
\[
S^1 \times [-1,1] \rightarrow \big( S^1 \times[-1,1] \big) \big/ \big( S^1 \times\{\pm1\} \big) \rightarrow M,
\]
where the first map is the natural quotient map and the second map is an embedding.
We denote by $\bar{o}$ the image of $\mathfrak{p}$, i.e., the pinched torus, and by $\{p_{0}\}=\mathfrak{p}\big(S^{1}\times\{\pm1\}\big)$ the pinch point.

Let $B$ be an open contractible neighbourhood of $p_0$ in $M$.
By the continuity of $\mathfrak{p}$, there is $a\in(0,1)$ sufficiently close to $1$, such that $L = \mathfrak{p}(S^1 \times [a,1]) \cup \mathfrak{p}(S^1 \times [-1,-a])$ is contained in $B$.
Consider the loops $\lambda_\pm$ on the Lagrangian submanifold $\bar{o}$ obtained as the composition $S^1 \hookrightarrow S^1 \times \{\pm a\} \xrightarrow{\mathfrak{p}} \lambda_\pm$.
Since $B$ is contractible, $\Gamma_J^2|_B$ is a trivial $S^1$ bundle, $\Gamma_J^2|_B \simeq B \times S^1$, and there is a section $\mathfrak{s}: B \to \Gamma_J^2|_B$.
Denote by $\mu_\pm$ the Maslov indices of $\lambda_\pm$, see Definition~\ref{def:maslov-index-lagrangian}.
Moreover, notice that since $B$ is contractible the Maslov indices $\mu_\pm$ do not depend on the choice of section, see Proposition~\ref{prop:maslov-index-simply-connected}.

\begin{thm}
\label{thm:.Maslov-bundle trivializable over pinched-torus}
The bundle $\Gamma_{J}^{2}|_{\bar{o}}$ is isomorphic to $\bar{o}\times S^1$ if and only if $\mu_+=\mu_-$.
\end{thm}

\begin{proof}
Since $\bar{o}$ is a subspace of $M$ and the bundle $\Gamma_J^2$ is locally trivializable over $M$, $\Gamma_J^2|_{\bar{o}}$ is locally trivializable over $\bar{o}$. 
To show that $\Gamma_J^2|_{\bar{o}} \simeq \bar{o} \times S^1$ we need to construct a global section on $\bar{o}$ by extending the local section $\mathfrak{s}|_{\bar{o} \cap B} : \bar{o} \cap B \to \Gamma_J^2|_{\bar{o} \cap B}$. 
We resort to another section on $o = \bar{o} \setminus \{p_0\}$. 
At each point $q \in o$, the tangent space $T_q o$ is a Lagrangian space.
Therefore, the assignment $q \mapsto T_q o$ gives a section $\tilde{\tau}: o \to \Lambda_J|_o$, which then induces a section $\tau: o \to \Gamma_J^2|_o$. 

There is a continuous map $h_-: \lambda_- \to S^1$ such that $\mathfrak{s}(q) = \tau(q) \cdot h_-(q)$ for all $q \in \lambda_-$, where the operation on the right-hand side is the inherent $S^1$ action on the principal bundle $\Gamma_J^2$. Similarly, there is a continuous map $h_+: \lambda_+ \to S^1$ such that $\mathfrak{s}(q) = \tau(q) \cdot h_+(q)$ for all $q \in \lambda_+$. Then $\deg h_\pm = - \mu_\pm$.

Let $C = \mathfrak{p}(S^1 \times [-a,a])$ with $\partial C = \lambda_- \sqcup \lambda_+$. Since $\lambda_-$ is a deformation retract of $C$, the map $h_-$ can be extended continuously to a map $h: C \to S^1$ such that $h|_{\lambda_-} = h_-$. Moreover, the extended map $h$ can chosen to satisfy $h|_{\lambda_+} = h_+$ if and only if $\deg h_+ = \deg h_-$, that is, if and only if $\mu_+ = \mu_-$.

Recall that $L = \mathfrak{p}(S^1 \times [a,1]) \cup \mathfrak{p}(S^1 \times [-1,-a])$ and notice that $L \cap C = \lambda_- \sqcup \lambda_+$.
Therefore, if $\mu_+ = \mu_-$ we can define the continuous section $\sigma: \bar o \to \Gamma_J^2|_{\bar o}$ by $\sigma(q) = \mathfrak{s}(q)$ for $q \in L$ and $\sigma(q) = \tau(q) \cdot h(q)$ for $q \in C$.
Conversely, if there is a continuous section $\sigma : \bar o \to \Gamma_J^2|_{\bar o}$, it follows that $\deg h_+ = \deg h_-$ and thus $\mu_+ = \mu_-$.
\end{proof}

Given an integrable Hamiltonian system $(M,\Phi,F)$, a compatible Hamiltonian vector field $X$, is a Hamiltonian vector field for which the flow $\varphi_{X}$ of $X$ commutes with the $\mathbb{R}^{2}$-action $\Phi$ and preserves the integral map $F$. 
That is, $F\circ\varphi_{X}^{t}=F$ and $\varphi_{X}^{t} \circ \Phi^{u} = \Phi^{u} \circ \varphi_{X}^{t}$ for $t\in\mathbb{R}$ and $u\in\mathbb{R}^{2}$. 
Since $o\cong S^{1}\times\mathbb{R}$, a compatible Hamiltonian vector field $X_{\theta}$ can be defined on a neighbourhood of $\bar{o}$ such that its flow $\varphi_{\theta}$ has period $1$ in a neighbourhood of $o$.
Exploiting $\varphi_{\theta}$ we can prove the following theorem about the triviality of the restricted bundle $\Gamma_{J}^{2}\big|_{\bar{o}}$.

\begin{thm}
\label{thm: Maslov S^1 Bundle over pinched tori in integrable systems}
If $\bar{o} = o\cup\{p_0\}$ is a singular orbit of an integrable Hamiltonian system $(M,\Phi,F)$, then $\mu_+=\mu_-$ and thus the bundle $\Gamma_J^2|_{\bar{o}}$ is isomorphic to $\bar{o}\times S^1$.
\end{thm}

\begin{proof}
Since $\varphi_{\theta}$ is periodic in a neighbourhood of $o$, the tangent maps $\varphi_{\theta,*}^t$ give a closed path 
\begin{equation}
\label{eq:loops of LG planes}
\mathrm{P}_w : [0,1] \ni t \mapsto \varphi_{\theta,*}^t(w) \in \Lambda_J
\end{equation}
in $\Lambda_J$ for any point $p \in o$ and $w \in \Lambda_J \big|_p$.
Since $p_0$ is a limit point of $o$, points on the fiber $\Lambda_J|_{p_0}$ are also limit points of the set $\Lambda_J|_o$. 
Since $\Lambda_J$ is Hausdorff, the mapping in \eqref{eq:loops of LG planes} also gives closed paths for each $w \in \Lambda_J|_{p_0}$.

Let $p_-$ and $p_+$ be points in $\lambda_-$ and $\lambda_+$, respectively.
Since $\Lambda_J|_{p_-} \cong \mathbb{U}(n) / \mathbb{O}(n)$ is path-connected, for any Lagrangian planes $w_-$ and $w'_-$, the loops
\[
t \mapsto det_{\mathbb{C}}^2 \circ \varphi_{\theta,*}^t(w_-)
\]
and
\[
t \mapsto det_{\mathbb{C}}^{2} \circ \varphi_{\theta,*}^t(w'_-)
\]
have the same degree. 
Since $o$ is invariant under the flow $\varphi_\theta^t$, we can particularly choose $w_-$ to be tangent to $o$ and see that the degree equals $\mu_-$. 
Therefore, for any $w_- \in \Lambda_J|_{p_-}$, the degree of the loop $t \mapsto det_{\mathbb{C}}^2 \circ \varphi_{\theta,*}^{t}(w_-)$ is $\mu_-$. 
The same conclusion can be drawn for $\mu_+$ and any $w_+ \in \Lambda_J|_{p_+}$.

We connect $p_-$ and $p_+$ with a path $\gamma$ in $\bar{o} \cap B$ such that $\gamma(-1) = p_-$, $\gamma(0) = p_0$ and $\gamma(1) = p_+$.
Since $\Lambda_J|_\gamma \cong \gamma \times \mathbb{U}(n)/\mathbb{O}(n)$, $\gamma$ can be lifted to a path $\tilde{\gamma}$ in $\Lambda_J|_\gamma$ with $\tilde{\gamma}(-1) = w_-$ and $\tilde{\gamma}(+1) = w_+$.
Then 
\[
h_s(t) := det_{\mathbb{C}}^2 \circ \varphi_{\theta,*}^t (\tilde{\gamma}(s))
\]
gives a continuous family of closed loops in $S^1$ indexed by $s \in [-1,1]$,
from which we deduce that $\mu_+ = \mu_-$.
\end{proof}

\begin{rem}
Notice that the equality $\mu_+ = \mu_-$ also follows from the fact that these are the Maslov indices of the $S^1$ action $\varphi_\theta$ and thus they are both equal to the local Maslov index of the $S^1$ action at the fixed point $p_0$.
\end{rem}

\section{Simultaneous axial rotations on \texorpdfstring{$S^2 \times S^2$}{S2xS2}}
\label{sec:S^2 times S^2}
\label{sec:S2xS2}

In Section~\ref{sec:example-S2}, we considered $S^2$ as a symplectic manifold with symplectic form $\omega_{S^2}$, Eq.~\eqref{eq:omega_S2}, and we considered the corresponding Maslov circle bundle $\Gamma_{S^2}$. 
The space $S^{2} \times S^{2}$ has the natural symplectic structure
$\omega = \omega_{S^2} \oplus \omega_{S^2}$, where for $\hat u = (u_1,u_2), \hat v =(v_1,v_2) \in T_{(p_1,p_2)} (S^2 \times S^2) = T_{p_1} S^2 \times T_{p_2} S^2$, we have
\[
\omega_{(p_1,p_2)}(\hat u, \hat v) 
= (\omega_{S^2} \oplus \omega_{S^2})_{(p_1,p_2)}((u_1,u_2),(v_1,v_2))
= p_1 \cdot (u_1 \times v_1) + p_2 \cdot (u_2 \times v_2).
\]
Moreover, we have the corresponding compatible Riemannian structure $g$
\[
g_{(p_1,p_2)}(\hat u, \hat v) 
= (g_{S^2} \oplus g_{S^2})_{(p_1,p_2)}((u_1,u_2),(v_1,v_2))
= u_1 \cdot v_1 + u_2 \cdot v_2,
\]
and almost complex structure
\[
J_{(p_1,p_2)}(\hat u) 
= (J_{S^2})_{p_1}(u_1) + (J_{S^2})_{p_2}(u_2)
= - p_1 \times u_1 - p_2 \times u_2.
\]
The unitary frame bundle $Fr_{S^2 \times S^2}^u$ can be described---using the set-theoretic definition adopted in this paper---as the set
\begin{align*}
Fr_{S^{2}\times S^{2}}^{u}
= \bigl\{ & (p_1, p_2) \in S^2 \times S^2,\,
     (\hat u_1, \hat u_2, \hat v_1, \hat v_2)
     \in (T_{p_1}S^2 \times T_{p_2}S^2)^4; \\
     & \;
     \hat v_i = -J \hat u_i;\;
     g(\hat u_i, \hat u_j) = \delta_{ij};\;
     g(\hat u_i, \hat v_j) = 0
     \bigr\},
\end{align*}
and then the corresponding Maslov $S^1$ bundle is $\Gamma_{S^2\times S^2} = Fr_{S^2 \times S^2}^u \big/ \mathbb{SU}(2)$.

In this section, we first construct $\Gamma_{S^2\times S^2}$ from $\Gamma_{S^2} \times \Gamma_{S^2}$, and then we consider the Maslov data for a specific $S^1$ action that is motivated by physical problems such as the perturbed Kepler problem and coupled angular momenta.

\subsection{The construction of \texorpdfstring{$\Gamma_{S^2 \times S^2}$}{Gamma(S2xS2)} from \texorpdfstring{$\Gamma_{S^2} \times \Gamma_{S^2}$}{Gamma(S2) x Gamma(S2)}}

Recall from Section~\ref{sec:example-S2} that $\Gamma_{S^2}$ is identified with the unitary frame bundle over $S^2$, and thus each point in $\Gamma_{S^2}$ can be seen as an orthonormal frame $w=(u,v)$ of the tangent space $T_p S^2$ for some $p\in S^2$.
The inherent $S^1$ action on the principal circle bundle $\Gamma_{S^2}$ is given by Eq.~\eqref{eq:gamma-s2-inherent-s1-action}.
Each point in $\Gamma_{S^2}\times\Gamma_{S^2}$ takes the form $(w_1,w_2)=(u_1,v_1,u_2,v_2)$ with $w_i = (u_i,v_i)\in\Gamma_{S^2}$ being an orthonormal frame of the tangent space $T_{p_{i}} S^{2}$ for some $p_i\in S^2$.
Moreover, the space $\Gamma_{S^2}\times\Gamma_{S^2}$ carries a natural $\mathbb T^2$ action given by
\[
(w_1,w_2) \cdot (e^{i t_1}, e^{i t_2}) 
= (w_1 \cdot e^{i t_1}, w_2 \cdot e^{i t_2}),
\]
where the $S^1$ actions at the right-hand side are given by Eq.~\eqref{eq:gamma-s2-inherent-s1-action}.
The $\mathbb{T}^2$ action on $\Gamma_{S^2}\times\Gamma_{S^2}$ can also be written as
\begin{equation}
(u_1, v_1, u_2, v_2) \cdot (e^{i t_1}, e^{i t_2}) 
= (u_1, v_1, u_2, v_2) \begin{bmatrix}
\cos t_1 & -\sin t_1 & 0 & 0\\
\sin t_1 & \cos t_1 & 0 & 0\\
0 & 0 & \cos t_2 & -\sin t_2\\
0 & 0 & \sin t_2 & \cos t_2
\end{bmatrix}
\label{eq:T2-action-GS2xGS2}
\end{equation}

\begin{prop}
The Maslov $S^1$ bundle $\Gamma_{S^2 \times S^2}$ is isomorphic to the circle bundle $(\Gamma_{S^2} \times \Gamma_{S^2}) \big/ \psi_{S^1}$,
where $\psi_{S^1}$ denotes the $S^1$ action 
\begin{equation}
\label{eq:psi_S^1}
\psi_{S^1}: S^1 \times \Gamma_{S^2}\times\Gamma_{S^2} \to \Gamma_{S^2}\times\Gamma_{S^2} 
: e^{it}, (w_1,w_2) \mapsto \psi_{S^1}^t(w_1,w_2)
= (w_1 \cdot e^{it}, w_2 \cdot e^{-it}).
\end{equation}

\end{prop}

\begin{proof}
The map $\iota : \Gamma_{S^2} \times \Gamma_{S^2} \to Fr_{S^2\times S^2}^{u}$ defined by
\[
\iota|_{(p_1,p_2)}(w_1 = (u_1,v_1), w_2 = (u_2,v_2)) 
= (\hat u_1 = (u_1,0), \hat u_2 = (0,u_2), \hat v_1 = (v_1,0), \hat v_2 = (0,v_2)),
\]
is an embedding of $\Gamma_{S^2} \times \Gamma_{S^2}$ into $Fr_{S^2\times S^2}^u$.
The following composition of maps, where $q_{\mathbb{SU}(2)}$ is the quotient map $Fr_{S^2 \times S^2}^u \to Fr_{S^2 \times S^2}^u / \mathbb{SU}(2)$, is then natural:
\begin{equation}
\kappa : \Gamma_{S^2} \times \Gamma_{S^2} 
\xrightarrow{\iota} Fr_{S^2 \times S^2}^u
\xrightarrow{q_{\mathbb{SU}(2)}} \Gamma_{S^2 \times S^2}.
\label{eq:kappa}
\end{equation}
The $\mathbb{T}^2$ action on $\Gamma_{S^2} \times \Gamma_{S^2}$ induces through the map $\iota$ a $\mathbb{T}^2$ action on $\iota(\Gamma_{S^2} \times \Gamma_{S^2})$ given by
\begin{equation}
\iota(w_1,w_2) \cdot (e^{i t_1}, e^{i t_2})
:= \iota(w_1 \cdot e^{i t_1}, w_2 \cdot e^{i t_2})
= \iota(w_1,w_2) \widehat R( t_1, t_2),
\label{eq:T^2 action on Fr^u}
\end{equation}
where
\[ \widehat R( t_1, t_2) = \begin{bmatrix}
\cos t_1 & 0 & -\sin t_1 & 0\\
0 & \cos t_2 & 0 & -\sin t_2\\
\sin t_1 & 0 & \cos t_1 & 0\\
0 & \sin t_2 & 0 & \cos t_2
\end{bmatrix} \in \mathbb{U}(2). 
\]
The map $\iota$ is $\mathbb{T}^2$-equivariant
with $\mathbb{T}^{2}$ acting on $Fr_{S^{2}\times S^{2}}^{u}$ as
a subgroup of $\mathbb{U}(2)$ via the embedding $(e^{i t_1},e^{i t_2}) \mapsto \widehat R( t_1, t_2)$. 
We check that
\begin{equation}
\begin{aligned}
\kappa(w_1\cdot e^{i t_1}, w_2\cdot e^{i t_2})
& = \kappa((w_1, w_2) \cdot (e^{i t_1}, e^{i t_2})) \\
& = q_{\mathbb{SU}(2)}(\iota((w_1, w_2) \cdot (e^{i t_1}, e^{i t_2}))) \\
& = q_{\mathbb{SU}(2)} (\iota(w_1,w_2) \widehat R( t_1, t_2)) \\
& = q_{\mathbb{SU}(2)}(\iota(w_1,w_2)) {\textstyle\det_{\mathbb C}} \widehat R( t_1, t_2) \\
& = \kappa(w_1,w_2) e^{i( t_1+ t_2)}.
\end{aligned}
\label{eq:extended-S1-equivariance-kappa}
\end{equation}

Since $\kappa$ is a bundle map that covers the identity map on the base space $S^2 \times S^2$, we have $\kappa(w_1',w_2') = \kappa(w_1,w_2)$ only if $(w_1',w_2') = (w_1,w_2) \cdot (e^{it_1}, e^{it_2})$. 
Then, Eq.~\eqref{eq:extended-S1-equivariance-kappa} implies that $\kappa(w_1',w_2') = \kappa(w_1,w_2) e^{i(t_1+t_2)}$ and we conclude that $t_2=-t_1$. 
Therefore, $\kappa(w_1',w_2') = \kappa(w_1,w_2)$ if and only if $(w_1',w_2') = \psi_{S^1}^t(w_1,w_2) = (w_1,w_2) \cdot (e^{it},e^{-it})$, that is, if and only if $[w_1',w_2'] = [w_1,w_2]$, where $[w_1,w_2]$ denotes the $\psi_{S^1}$ orbit through $(w_1,w_2)$.
This implies that $\kappa$ factorizes as
\begin{equation}
\kappa : \Gamma_{S^2} \times \Gamma_{S^2} 
\xrightarrow{q_{S^1}} (\Gamma_{S^2} \times \Gamma_{S^2}) \big/ \psi_{S^1}
\xrightarrow{\mathcal{I}} \Gamma_{S^2 \times S^2},
\label{eq:factorization-kappa}
\end{equation}
where $q_{S^1}$ is the quotient map for the $S^1$ bundle $\Gamma_{S^2} \times \Gamma_{S^2} \to (\Gamma_{S^2} \times \Gamma_{S^2}) \big/ \psi_{S^1}$, and $\mathcal{I}$ satisfies $\mathcal{I}([w_1,w_2]) = \kappa(w_1,w_2)$ for any representative $(w_1,w_2)$ of $[w_1,w_2]$.
In particular, we have the commutative diagram
\begin{equation*}
\begin{tikzcd}
\Gamma_{S^2}\times\Gamma_{S^2} \rar["\iota"] \drar["\kappa"] \dar["q_{S^1}"] & Fr^u_{S^2 \times S^2} \dar["q_{\mathbb{SU}(2)}"] \\
(\Gamma_{S^2}\times\Gamma_{S^2}) / \psi_{S^1} \rar["\mathcal{I}"] & \Gamma_{S^2 \times S^2}.
\end{tikzcd}
\end{equation*}

We now show that $\mathcal{I}$ is an isomorphism between the principal $S^1$
bundles $\big(\Gamma_{S^2}\times\Gamma_{S^2}\big)\big/\psi_{S^1}$
and $\Gamma_{S^2\times S^2}$.
The map $\mathcal I$ is injective, since if $\mathcal{I}([w_1',w_2']) = \mathcal{I}([w_1,w_2])$, then $\kappa(w_1',w_2') = \kappa(w_1,w_2)$ and thus $[w_1',w_2'] = [w_1,w_2]$.
Moreover, since $\kappa$ is a bundle morphism that covers the identity map on the base space $S^{2}\times S^{2}$, Eq.~\eqref{eq:extended-S1-equivariance-kappa} 
implies that $\kappa$ is a surjection, and then so is $\mathcal{I}$.

The space $\big(\Gamma_{S^{2}}\times\Gamma_{S^{2}}\big)\big/\psi_{S^{1}}$ is a principal $S^1$ bundle with inherent $S^1$ action given by
\[
[w_{1},w_{2}]\cdot e^{it}=[w_{1}\cdot e^{it},w_{2}]=[w_{1},w_{2}\cdot e^{it}].
\]
Finally, using ~Eq.~\eqref{eq:extended-S1-equivariance-kappa} we check that
\[
\mathcal{I}([w_1,w_2]\cdot e^{it})
= \mathcal{I}([w_1 \cdot e^{it}, w_2])
= \kappa(w_1 \cdot e^{it},w_2) 
= \kappa(w_1, w_2) e^{it}
= \mathcal{I}([w_1, w_2]) e^{it}
,
\]
and hence $\mathcal{I}$ is a bundle isomorphism.
\qedhere

\end{proof}

To define a connection $1$-form for the principal bundle $\Gamma_{S^2 \times S^2}$, 
we first note that, the vertical distribution $\mathcal{V}(\Gamma_{S^2}\times\Gamma_{S^2})$ is naturally isomorphic to the direct sum $\mathcal{V}\Gamma_{S^2}\oplus\mathcal{V}\Gamma_{S^2}$,
and, the tangent bundle $T(\Gamma_{S^2}\times\Gamma_{S^2})$ is isomorphic to $T\Gamma_{S^2}\oplus T\Gamma_{S^2}$. 
Denote by $\partial/\partial\theta_1$ and $\partial/\partial\theta_2$ the infinitesimal generators of the inherent $\mathbb T^2$ action on $\Gamma_{S^2} \times \Gamma_{S^2}$ corresponding to each of the factors $\Gamma_{S^2}$, 
and by $\partial/\partial\theta$ the infinitesimal generator of the inherent $S^1$ action on $\Gamma_{S^2 \times S^2}$.

Let $\beta=f_{\beta}\cdot\frac{\partial}{\partial\theta}$ be a connection $1$-form on $\Gamma_{S^2}$ with the horizontal distribution $\ker_{\beta}$. 
Since $\kappa$ is a bundle morphism that covers the identity map on the base $S^{2}\times S^{2}$, the restriction of the tangent map $\kappa_*$ to the subspace $\ker_{\beta}\oplus\ker_{\beta}$ is non-degenerate. 
Eq.~\eqref{eq:extended-S1-equivariance-kappa} implies that 
\begin{equation}\label{eq:kappa-tangent}
\kappa_* \Bigl( \lambda_1 \frac{\partial}{\partial\theta_1},  \lambda_2 \frac{\partial}{\partial\theta_2} \Bigr) 
= (\lambda_1+\lambda_2) \frac{\partial}{\partial\theta}.
\end{equation}
Therefore, $\ker\kappa_*$ is exactly the $1$-dimensional vector bundle spanned by the infinitesimal generator $X_{\psi} = (\partial/\partial\theta_1,-\partial/\partial\theta_2)$ of the $S^1$ action $\psi_{S^1}$. 
Define the $1$-form $f_\beta \oplus f_\beta$ on $\Gamma_{S^2} \times \Gamma_{S^2}$
by
\begin{equation*}
f_{\beta}\oplus f_{\beta}(u,v)=f_{\beta}(u)+f_{\beta}(v),
\end{equation*}
for $(u,v)\in T_{w_1}\Gamma_{S^{2}}\oplus T_{w_2}\Gamma_{S^{2}}$.

\begin{prop}
There is a unique $1$-form $f_{\beta^\oplus}$ on $\Gamma_{S^2 \times S^2}$ satisfying 
\begin{equation}
f_\beta \oplus f_\beta = \kappa^* f_{\beta^\oplus}.
\label{eq:f-beta-oplus}
\end{equation}
The $1$-form $\beta^\oplus = f_{\beta^\oplus} \cdot \frac{\partial}{\partial\theta}$ is a connection $1$-form on $\Gamma_{S^2 \times S^2}$.
\end{prop}

\begin{proof}
It is straightforward to check that $\ker\kappa_{*}\subset\ker f_{\beta} \oplus f_{\beta}$,
and then due to the surjectivity of $\kappa_{*}$, 
there exists a unique $1$- form $f_{\beta^\oplus}$ on $\Gamma_{S^2 \times S^2}$
such that
\begin{equation*}
f_\beta \oplus f_\beta = f_{\beta^\oplus} \circ \kappa_*,
\end{equation*}
which is Eq.~\eqref{eq:f-beta-oplus}.

Eq.~\eqref{eq:kappa-tangent} gives $\kappa_*(\partial/\partial\theta_1,0) = \partial/\partial\theta$,
and then as a result,
\[
f_{\beta^{\oplus}}\left(\frac{\partial}{\partial\theta}\right)
= f_{\beta}\left(\frac{\partial}{\partial\theta_1}\right)
= 1.
\]
Since the restriction $\kappa_*|_{\ker_{\beta}\oplus\ker_{\beta}}$
is nondegenerate, dimensionality implies that
\[
\kappa_*(\ker_\beta \oplus \ker_\beta)=\ker f_{\beta^\oplus},
\]
and Eq.~\eqref{eq:extended-S1-equivariance-kappa} implies the invariance of $\ker f_{\beta^{\oplus}}$ under the inherent $S^1$ action on $\Gamma_{S^{2}\times S^{2}}$. As a result, 
\[
\beta^{\oplus}=f_{\beta^{\oplus}}\cdot\frac{\partial}{\partial\theta}
\]
defines a connection $1$-form on $\Gamma_{S^{2}\times S^{2}}$.
\end{proof}

\subsection{Local Maslov indices of a symplectic \texorpdfstring{$S^1$}{S1} action on \texorpdfstring{$S^2 \times S^2$}{S2xS2}}

In this section, we discuss the symplectic $S^1$ action on $S^2 \times S^2$ which simultaneously rotates the spheres about their ``$z-$axes'' and we compute the values of the local Maslov indices at the fixed points of the action.

To be explicit, viewing $S^2 \times S^2$ as a submanifold in $\mathbb R^3 \times \mathbb R^3$, we consider the $S^1$ action 
\begin{equation}
e^{it}, (p_1, p_2) \xmapsto{\varphi_\ax} 
(R_{z}^{kt} p_1, R_{z}^{lt} p_2),
\label{eq:axial-s1-action}
\end{equation}
where $p_1,p_2 \in S^2 \subseteq \mathbb R^3$, $R_{z}^{\theta}$ stands for the $3\times3$ matrix
\[
\begin{bmatrix}
\cos\theta & -\sin\theta & 0\\
\sin\theta & \cos\theta & 0\\
0 & 0 & 1
\end{bmatrix},
\]
and $k,l$ are coprime integers.

The action $\varphi_\ax$ is Hamiltonian and it has $4$ fixed points $(p_{\pm},p_{\pm})$, where $p_\pm = (0,0,\pm1)^{\mathrm{T}}$.
We compute the values of the local Maslov indices $\mathfrak{m}_{\pm,\pm}$ of $\varphi_\ax$ at the fixed points $(p_\pm,p_\pm)$.

The $S^1$ action $R_z$ on $S^2$ is lifted to an action $\mathcal{R}_{z}$
on the principal $S^{1}$ bundle $\Gamma_{S^2}$ with
\[
(u,v), e^{i\theta} \xmapsto{\mathcal{R}_{z}}(R_{z,*}^{\theta}u, R_{z,*}^{\theta}v),
\]
where $R_{z,*}^\theta$ is the tangent map of $R_z^\theta$.
Note that $\mathcal{R}_{z}$ commutes with the inherent $S^1$ action
on $\Gamma_{S^{2}}$, and hence on the fibers $\Gamma_{S^{2}}|_{p_{\pm}}$
over the fixed points $p_{\pm}$ of $R_{z}$, there exist constants
$c_{\pm}\in\mathbb{R}$ such that
\[
\mathcal{R}_{z}^{e^{i\theta}} w = w\cdot e^{ic_{\pm}\theta}
\]
for $w = (u,v) \in \Gamma_{S^2}|_{p_\pm}$. 
The orbits $\mathcal{R}_{z}^{\theta}w$ are closed and thus they wind around $\Gamma_{S^2}|_{p_{\pm}}\cong S^{1}$ an integer number of times, implying that $c_{\pm}$ can only be integers.
One then easily checks that $c_\pm = \pm1$.

Now we look into the symplectic $S^1$ action $\varphi_\ax$. 
Recall from Sec.~\ref{sec:compact-group-actions}, that $\varphi_\ax$ is lifted to an action $\Phi_\ax$ on $Fr_{S^{2}\times S^{2}}^{u}$ given by
\[
(\hat u_1, \hat u_2, \hat v_1, \hat v_2), e^{it}
\xmapsto{\Phi_\ax}
\big(
\varphi_{\ax,*}^{e^{it}} (\hat u_1),
\varphi_{\ax,*}^{e^{it}} (\hat u_2), 
\varphi_{\ax,*}^{e^{it}} (\hat v_1),
\varphi_{\ax,*}^{e^{it}} (\hat v_2)
\big).
\]

Note that when viewing the tangent bundle of $S^2 \times S^2$ as $TS^2 \times TS^2$, the tangent map $\varphi_{\ax,*}^{e^{it}}$
is simply the map sending $\hat u = (u_1,u_2) \in TS^2 \times TS^2$ to $(R_{z,*}^{kt} u_1, R_{z,*}^{lt}u_2)$.
It is then straightforward to check that 
\[
\Phi_\ax^{e^{it}} \circ \iota = \iota \circ 
(\mathcal{R}_{z}^{e^{ikt}} \times \mathcal{R}_{z}^{e^{ilt}}). 
\]
We denote the lifted action of $\varphi_\ax$ on $\Gamma_{S^{2} \times S^{2}}$ also by $\Phi_\ax$, and we then have the following commutative diagram:
\[
\begin{CD}
\Gamma_{S^{2}}\times\Gamma_{S^{2}} & @>\iota>> & Fr_{S^{2}\times S^{2}}^{u} & @> q_{\mathbb{SU}(2)} >> & \Gamma_{S^{2}\times S^{2}}\\
@V\mathcal{R}_{z}^{e^{ikt}}\times\mathcal{R}_{z}^{e^{ilt}}VV &  & @V\Phi_{\ax}^{e^{it}}VV &  & @V\Phi_{\ax}^{e^{it}}VV\\
\Gamma_{S^{2}}\times\Gamma_{S^{2}} & @>\iota>> & Fr_{S^{2}\times S^{2}}^{u} & @> q_{\mathbb{SU}(2)} >> & \Gamma_{S^{2}\times S^{2}}
\end{CD}.
\]
This implies $\Phi_{\ax}^{e^{it}} \circ \kappa = \kappa \circ (\mathcal{R}_{z}^{e^{ikt}}\times\mathcal{R}_{z}^{e^{ilt}})$. 
Therefore, for $a,b\in\{+,-\}$, we have the following relation over the fixed point $(p_{a},p_{b})$
\[
\Phi_{\ax}^{e^{i\theta}}(\kappa(w_1,w_2))
= \kappa(\mathcal{R}_z^{e^{ik\theta}} w_1, \mathcal{R}_z^{e^{il\theta}} w_2)
= \kappa(w_1 \cdot e^{ic_{a}k\theta}, w_2 \cdot e^{ic_{b}l\theta}).
\]
From Eq.~\eqref{eq:extended-S1-equivariance-kappa} we get
\[
\Phi_{\ax}^{e^{i\theta}} (\kappa(w_1,w_2))
= \kappa(w_1,w_2) e^{i(c_{a}k+c_{b}l)\theta}.
\]

The local Maslov index $\mathfrak{m}_{a,b}$ of $\varphi_{\ax}$ at
$(p_{a},p_{b})$ is twice the number of times the orbit of $\Phi_{\ax}$ winds around the fiber $\Gamma_{(p_{a},p_{b})}$, and hence
\[
\mathfrak{m}_{a,b} = 2 (c_a k + c_b l),
\]
that is, $\mathfrak{m}_{+,+} = -\mathfrak{m}_{-,-} = 2(k+l)$, and
$\mathfrak{m}_{+,-} = - \mathfrak{m}_{-,+} = 2(k-l)$.

\section{Conclusions}
\label{sec:conclusions}

In this work, we defined the principal $S^1$ bundles $\Gamma_J$ and $\Gamma_J^2$ which we call Maslov bundles.
When the Maslov bundles are trivializable, they are directly related to the usual notion of Maslov indices.
We used these bundle structures to analyze the symplectic actions of compact groups on symplectic manifolds via the lifted actions on the bundles. 

In particular, we investigated symplectic actions on homogeneous spaces and monotone manifolds.
For homogeneous spaces, we showed that the first real Chern class being nonzero implies $\Gamma_J^2$ and $\Gamma_J$ to be homogeneous spaces as well. 
For monotone manifolds, we used Maslov $S^1$ bundles to revisit Ono's result about symplectic actions being Hamiltonian \cite{Ono1988}. 
It is worth pointing out that, while \cite{Ono1988} demonstrated the result for an $S^1$ action on a compact monotone manifold, as is shown in this paper, the result also applies to any symplectic compact group action on a general (not necessarily compact) monotone symplectic manifold.

Furthermore, we defined the notion of the local Maslov index of a symplectic $S^1$ action and its relation to the usual Maslov index when $\Gamma_J^2$ is trivializable. 
We have shown that the local Maslov indices of a symplectic $S^1$ action are equal to the values of the corresponding momentum at the fixed points.
This directly implies that, on a compact monotone symplectic manifold with zero Chern class, the only symplectic $S^1$ action is the trivial action with constant momentum.

Besides the results mentioned above, we discussed the impact
of local Maslov indices of an $S^1$ action on the Hamiltonian monodromy of an integrable Hamiltonian system.
In particular, for a two degrees of freedom system with a compatible
$S^1$ action, a nonzero local Maslov index at a fixed point of the $S^1$ action implies the existence of a second $S^1$ action on the phase space (excluding possibly the singular leaves), and thus also implies the triviality of Hamiltonian monodromy.

Finally, we gave a detailed discussion of the structure of the Maslov $S^1$ bundles over Lagrangian pinched tori and over $S^{2}\times S^{2}$. 
These two spaces are of particular interest in the study of integrable Hamiltonian systems. 
The restriction of $\Gamma_J^2$ over a Lagrangian pinched torus was shown to be trivializable when the pinched torus is a singular leaf of an integrable system. 
Additionally, the construction of $\Gamma_{S^{2}\times S^{2}}$ from $\Gamma_{S^{2}}\times\Gamma_{S^{2}}$ was described in detail, and the local Maslov indices of an $S^1$ action on $S^2 \times S^2$ were calculated. 

\appendix 

\section{Definitions and conventions}
\label{sec:notation}

We recall basic definitions and properties of principal $S^1$ bundles and unitary frames.

\subsection{Principal \texorpdfstring{$S^1$}{S1} bundles}
\label{sec:principal-s1-bundles}

A principal $S^1$-bundle consists of a total space $P$ on which there is a free $S^1$ action, its orbifold $B=P\big/S^{1}$, and the quotient map $\pi_{P}:P\rightarrow B$.
To distinguish it from other group actions on $P$, in the rest of the paper we will call this $S^1$ action the inherent $S^1$
action. 
For simplicity, we will also call $P$ the principal bundle
when the $S^1$ action, and the spaces $B$ and $P$ are clear from the context.

A connection $\mathcal{H}$ on the bundle $\pi_{P}:P\rightarrow B$ is a horizontal distribution invariant under the inherent $S^{1}$ action such that the tangent bundle $TP$ is split as a direct sum $\mathcal{H}\oplus\mathcal{V}$ with $\mathcal{V} = \ker(\pi_{P})_{*}$ being the vertical distribution.
The Lie algebra of $S^{1}$ is $T_{1}S^{1} = \mathbb{R}\cdot\frac{\partial}{\partial\theta}$
with $\frac{\partial}{\partial\theta}=\frac{d}{d\theta}\big|_{\theta=0}e^{i\cdot2\pi\theta}.$
The map
\[
\eta:P\times\Bigl(\mathbb{R}\cdot\frac{\partial}{\partial\theta}\Bigr) \to \mathcal{V},
\]
with 
\[
\Bigl( p, r\frac{\partial}{\partial\theta} \Bigr)
\mapsto \frac{d}{d\theta} \bigg|_{\theta=0}(p\cdot e^{i\cdot2\pi r\cdot\theta}),
\]
is a vector bundle isomorphism. 
For simplicity, we also denote the vector $\frac{d}{d\theta}\big|_{\theta=0}(p\cdot e^{i\cdot2\pi r\cdot\theta})$
by $r\frac{\partial}{\partial\theta}$. 
The connection $1$-form associated to the connection $\mathcal{H}$ is
\[
\alpha : TP \xrightarrow{pr_{\mathcal{V}}} \mathcal{V} 
\cong P\times\Bigl(\mathbb{R}\cdot\frac{\partial}{\partial\theta}\Bigr) 
\to \mathbb{R} \cdot \frac{\partial}{\partial\theta}
\]
with $pr_{\mathcal{V}}$ being the projection associated to the splitting $TP=\mathcal{H}\oplus\mathcal{V}$ onto $\mathcal{V}$.
In particular, $\alpha(\frac{\partial}{\partial\theta}) = \frac{\partial}{\partial\theta}$ and $\ker\alpha = \mathcal H$.

The connection $1$-form $\alpha$ is $S^{1}$ invariant, and there is an $S^{1}$-invariant $1$-form $f_{\alpha}$ such that $\alpha = f_{\alpha} \cdot \frac{\partial}{\partial\theta}$.
Note that $f_{\alpha}(\frac{\partial}{\partial\theta})\equiv1$ and $\ker f_{\alpha}=\mathcal{H}$.
With a slight abuse in terminology, we also call $f_{\alpha}$ a connection $1$-form on $P$.

The $S^{1}$-invariance of $f_{\alpha}$ implies $\mathcal{L}_{\frac{\partial}{\partial\theta}}f_{\alpha}=0$.
Then Cartan's formula gives
\[
0 = \mathcal{L}_{\frac{\partial}{\partial\theta}} f_{\alpha} = \iota_{\frac{\partial}{\partial\theta}} df_{\alpha} + d\Bigl(f_{\alpha}\Bigl(\frac{\partial}{\partial\theta}\Bigr)\Bigr) = \iota_{\frac{\partial}{\partial\theta}} df_{\alpha}.
\]
The last equation implies that 
\[
df_{\alpha}(pr_{\mathcal{H}}\cdot,pr_{\mathcal{H}}\cdot)=df_{\alpha},
\]
where $pr_{\mathcal{H}}$ is the projection of $\mathcal{H}\oplus\mathcal{V}$ onto $\mathcal{H}$.
Since $df_{\alpha}$, $\mathcal{H}$, and $\mathcal{V}$, are all $S^{1}$-invariant, this implies that there is a closed $2-$form $\Omega_{\alpha}$ on $B$, called the curvature form, such that $\pi_{P}^{*}\Omega_{\alpha} = df_{\alpha}$.
Then, $[\Omega_{\alpha}] \in H_{dR}^{2}(B)$
is independent of $\alpha$, and is called the characteristic class
of the bundle $P \rightarrow B$.

We close this discussion with the following result which is used in the proof of Theorem~\ref{thm:homogeneous-G-space-(intro)}.

\begin{lem}
\label{lem:characteristic-class-integrable-connection}
The characteristic class $[\Omega_\alpha] \in H_{dR}^{2}(B)$ of a principal $S^1$ bundle vanishes if and only if the bundle admits an integrable connection.
\end{lem}

\begin{proof}
We first show that a connection $\mathcal{H} = \ker\alpha$ is integrable if and only if the corresponding curvature form $\Omega_\alpha = 0$.
According to the Frobenius integrability theorem, $\mathcal{H}$ is an integrable distribution if and only if $f_{\alpha} \wedge df_{\alpha}=0$.
The last equation is satisfied if only if $df_{\alpha} = 0$. 
This is because, if $df_{\alpha} \ne 0$, then there exist
$u,v\in\mathcal{H}_{p}=\ker f_{\alpha}\big|_{p}$ such that $df_{\alpha}(u,v)\neq0$,
and then 
\[
f_{\alpha}\wedge df_{\alpha}\Bigl(\frac{\partial}{\partial\theta},u,v \Bigr)
= f_{\alpha}\Bigl(\frac{\partial}{\partial\theta}\Bigr) \cdot df_{\alpha}(u,v) = df_{\alpha}(u,v) \ne 0,
\]
yielding a contradiction. Therefore, $\mathcal{H}$ is an integrable connection, if and only if  $df_\alpha = 0$, which from the surjectivity of $(\pi_P)_*$ is true if and only if $\Omega_{\alpha}=0$.

For the proof of the lemma, we have that if $\mathcal H$ is integrable then $\Omega_\alpha = 0$, and thus $[\Omega_\alpha] = 0$.
Conversely, if $[\Omega_\alpha] = 0$, there exists a $1$-form $\tau$ on $B$ such that $d\tau = \Omega_{\alpha}$. 
Then $f_{\alpha'} = f_\alpha - \pi_P^*(\tau)$ is a connection 1-form with curvature form $\Omega_{\alpha'} = 0$.
To see this, first notice that $f_{\alpha'}$ is $S^{1}$-invariant with
\[
f_{\alpha'}\Bigl(\frac{\partial}{\partial\theta}\Bigr)
= f_{\alpha}\Bigl(\frac{\partial}{\partial\theta}\Bigr) 
= 1,
\]
and 
\[
df_{\alpha'}
= \pi_P^*(\Omega_{\alpha}) - \pi_P^*(d\tau) = 0,
\]
implying that $\Omega_{\alpha'} = 0$ and that the corresponding distribution $\mathcal H' = \ker f_{\alpha'}$ is integrable.
\end{proof}

\subsection{Unitary frames}
\label{sec:unitary-frames}

We consider a symplectic manifold $(M, \omega)$ with $\dim M = 2n$ with an almost complex structure $J$ and a Riemannian structure $g$, so that $(\omega, J, g)$ are compatible, that is,%
\footnote{Our convention is different from the one used in \cite{CannasDaSilva2001} and corresponds to a sign change for the almost complex structure.}
\[ g(u, v)  = \omega(J u,v), \quad u, v \in T_pM. \]

A \emph{unitary frame} of $T_pM$ is a basis $(u_1, \dots, u_n, v_1, \dots, v_n)$ satisfying
\[ 
  v_i = - J u_i, \quad g(u_i,u_j) = \delta_{ij}, \quad g(u_i,v_j) = 0. 
\]
One can easily check that a unitary frame is both an orthonormal basis of $T_pM$ with respect to $g$, and a symplectic basis with respect to $\omega$.

Define the Hermitian inner product
\[ h(u,v) = g(u,v) + i \omega(u,v). \]
A linear transformation $U$ on $T_pM$ is \emph{unitary} with respect to $h$ if $h(Uu, Uv) = h(u, v)$.
This implies $g(Uu,Uv) = g(u,v)$ and $\omega(Uu,Uv) = \omega(u,v)$, that is, $U$ is orthogonal with respect to $g$ and symplectic with respect to $\omega$.
Moreover, one can show that $JU = UJ$ and thus $U$ also preserves the almost complex structure $J$.

A unitary linear transformation $U$ expressed in terms of a unitary frame $(u_1,\dots,u_n,v_1,\dots,v_n)$ has the matrix representation
\begin{equation}
\label{eq:unitary-matrix}
U = \begin{bmatrix}
A & -B \\ B & A
\end{bmatrix},
\end{equation}
where $A$, $B$ are real $n \times n$ matrices satisfying
\[ A^t A + B^t B = \operatorname{id}, \quad A^t B = B^t A. \]
Then the vectors $(u_1',\dots,u_n',v_1',\dots,v_n')$ defined by
\[ 
u_i' = U u_i = \sum_k A_{ki} u_k + B_{ki} v_k, \text{ and }
v_i' = U v_i = \sum_k -B_{ki} u_k + A_{ki} v_k,
\]
form another unitary frame of $T_pM$ with respect to $h$.

Write $w_i = u_i$, $w_{n+i} = v_i$ for $i=1,\dots,n$.
Denote by $W$ the matrix with the vector $w_i$ in its $i$-th column, that is, $W_{ji} = w_{i,j}$.
If $w'_i = U w_i$, $i=1, \dots, n$, then we have $w'_i = \sum_k U_{ki} w_k$ and thus
\[ W'_{ji} = w'_{i,j} = \sum_k U_{ki} w_{k,j} = \sum_k W_{jk} U_{ki} = (WU)_{ji}. \]
Therefore, $W' = WU$, that is, $\mathbb U(n)$ acts on the space of unitary frames on the right.

Given a unitary matrix $U$ as in Eq.~\eqref{eq:unitary-matrix}, the matrix $\widehat U = A + i B$ is a matrix in the classical complex unitary group $\mathbb U(n)$. 
The complex determinant of $U$ is defined by $\det_{\mathbb C} U = \det \widehat U = \det(A + i B)$.

\section{The topology of \texorpdfstring{$\Lambda_J$}{the bundle of Lagrangian planes}}
\label{sec:topology-lambda}

In this work we use the bundle of Lagrangian planes $\Lambda_{pl}$ interchangeably with the bundle $\Lambda_J = Fr_J^u / \mathbb{O}(n)$.
In this appendix we show that these spaces with their corresponding natural smooth structures are homeomorphic.

Consider the bundle $Fr^{\mathrm{sp}}$ of symplectic frames over $M$.
That is, over each $p\in M$, the fiber $Fr^{\mathrm{sp}}|_{p}$ consists of elements taking the form $(u_1, \dots, u_n, v_1, \dots, v_n)$ with $u_i,v_j\in T_pM$, such that $\omega(u_i,v_j)=\delta_{ij}$ and $\omega(u_i,u_j) = \omega(v_i,v_j) = 0$. 
Note that $Fr^{\mathrm{sp}}$ as well as $Fr^u_J$ have the natural bundle structures inherited from the frame bundle of $TM$ as its subbundles. 
The bundle $\Lambda_J$, being the quotient space $Fr_J^u / \mathbb{O}(n)$, inherits a topology from $Fr^u_J$.
For the bundle $\Lambda_{pl}$, consider the submersion $\Pi:Fr^{\mathrm{sp}} \to \Lambda_{pl}$ defined by
\[
\Pi(u_1, \dots, u_n, v_1, \dots, v_n) = \mathrm{span}\{u_1, \dots, u_n\},
\]
mapping each fiber $Fr^{\mathrm{sp}} |_p$ to the corresponding fiber $\Lambda_{pl} |_p$. 
Then $\Pi$ is surjective and hence it induces a quotient topology on $\Lambda_{pl}$.
To show that $\Lambda_{pl}$ and $\Lambda_J$ are homeomorphic we first construct a submersion $P_J : Fr^{\mathrm{sp}} \to Fr_J^u$.

We define the \emph{Gram-Schmidt map} $\mathcal G: \mathbb{SP}(2n) \to \mathbb{U}(n)$ by applying the following algorithm to $A \in \mathbb{SP}(2n)$. 
First, for $i=1,\dots,n$, define
\[ A'_i = A_i - \sum_{j=1}^{i-1} \frac{\langle A_i,A'_j\rangle}{\langle 
A'_j,A'_j\rangle} A'_j, \]
where $A_i$ is the $i$-th column of $A$ and $\langle A_i, A_j \rangle$ is the standard inner product in $\mathbb R^{2n}$.
Then, 
\[ \mathcal G(A) = \biggl[ 
\frac{A'_1}{|A'_1|}, \dots, \frac{A'_n}{|A'_n|},
-\frac{J_0 A'_1}{|A'_1|}, \dots, -\frac{J_0A'_n}{|A'_n|}
\biggr],
\]
where $J_0$ is the canonical almost complex structure on $\mathbb R^{2n}$ and $|A'_i|^2 = \langle A'_i, A'_i \rangle$.

\begin{lem}
The Gram-Schmidt map $\mathcal G: \mathbb{SP}(2n) \to \mathbb{U}(n)$ is left $\mathbb U(n)$ equivariant and as a result it is a submersion.
\end{lem}

\begin{proof}
A symplectic matrix $A \in \mathbb{SP}(2n)$ has a unique Iwasawa decomposition $A = K D N$ where $K \in \mathbb U(n)$, $D$ is a diagonal matrix with positive diagonal elements $d_1, \dots, d_n, d_1^{-1}, \dots, d_n^{-1}$, and $N$ is a nilpotent matrix having the form
\[ N = \begin{bmatrix} N_1 & N_2 \\ 0 & N_1^{-t} \end{bmatrix}, \]
where $N_1, N_2$ are $n \times n$ matrices satisfying $N_1 N_2^t = N_2 N_1^t$.
It is straightforward to check that $\mathcal G(A) = K$, that is, $A = \mathcal G(A) D N$.
If $C \in \mathbb U(n)$ and $A = \mathcal G(A) D N$, then $C A = C \mathcal G(A) D N$, and the uniqueness of the decomposition implies that $\mathcal G(CA) = C \mathcal G(A)$.
\end{proof}

For any open subset $U$ of $M$ over which the tangent bundle admits local trivializations, consider a smooth section $w : U \to Fr_J^u$ of unitary frames. The section $w$ induces local trivializations $\alpha_U: Fr_J^u|_U \to U \times \mathbb U(n)$ and $\beta_U: Fr^{\mathrm{sp}} |_U \to U \times \mathbb{SP}(2n)$.
Denote by $\mathcal G_U : U \times \mathbb{SP}(2n) \to U \times \mathbb{U}(n)$ the map defined by $\mathcal G_U(b, A) = (b, \mathcal G(A))$.
Then we can use the local maps $\mathcal G_U$ to construct a global submersion $\mathcal G_J$ from $Fr^{\mathrm{sp}}$ to $Fr_J^u$.

\begin{lem}
\label{lem:G_J_submersion}
The map $\mathcal{G}_J|_U : Fr^{\mathrm{sp}}|_U \to Fr_J^u |_U$ defined as $\mathcal G_J|_U = \alpha_U^{-1} \circ \mathcal G_U \circ \beta_U$ does not depend on the choice of trivializations and as a result it extends to a map $\mathcal{G}_J : Fr^{\mathrm{sp}} \to Fr_J^u$ which is a submersion.
\end{lem}

\begin{proof}
Consider another smooth section of unitary frames $w' : U \to Fr^u_J$ and corresponding trivializations $\alpha'_U: Fr_J^u|_U \to U \times \mathbb U(n)$ and $\beta'_U: Fr^{\mathrm{sp}} |_U \to U \times \mathbb{SP}(2n)$. For each $b \in U$, there is some $C \in \mathbb U(n)$ such that $w(b) = w'(b) C$.
For a symplectic frame $f^{\mathrm{sp}} \in Fr^{\mathrm{sp}}$ at $b \in U$ there is a unique $A \in \mathbb{SP}(2n)$ such that $f^{\mathrm{sp}} = w(b) A$. Then $\beta_U(f^{\mathrm{sp}}) = (b, A)$. Since $f^{\mathrm{sp}} = w'(b) C A$ we obtain $\beta'_U(f^{\mathrm{sp}}) = (b, C A)$. Similarly, if $f^u \in Fr^u_J$ is a unitary frame at $b \in U$ then there is a unique $K \in \mathbb U(n)$ such that $f^u = w(b) K = w'(b) C K$. Therefore, $\alpha_U(f^u) = (b, K)$ and $\alpha'_U(f^u) = (b, C K)$.

Let $f^u = \alpha_U^{-1} \circ \mathcal G_U \circ \beta_U (f^{\mathrm{sp}}) = \alpha_U^{-1}(b, \mathcal G(A))$,
where $f^{\mathrm{sp}} = w(b) A$. Then we compute
\[
\begin{aligned}
(\alpha'_U)^{-1} \circ \mathcal G_U \circ \beta'_U (f^{\mathrm{sp}})
& = (\alpha'_U)^{-1}(\mathcal G_U(b, C A))
= (\alpha'_U)^{-1}(b,\mathcal G(C A)) \\
& = (\alpha'_U)^{-1}(b, C \mathcal G(A)) 
= \alpha_U^{-1}(b, \mathcal G(A)) 
= f^u.
\end{aligned}
\]
Therefore, $\mathcal G_J|_U$ does not depend on the choice of trivialization and it defines a map $\mathcal G_J : Fr^{\mathrm{sp}} \to Fr^u$. 
Finally, $\mathcal G_J$ is a submersion since $\mathcal G$ is a submersion.
\end{proof}

\begin{rem}
It can be checked that the map $\mathcal G_J$ in Lemma~\ref{lem:G_J_submersion} corresponds to applying the following algorithm, paralleling the algorithm defining the Gram-Schmidt map $\mathcal G$, to a symplectic frame $f^{\mathrm{sp}} = (u_1,\dots,u_n,v_1,\dots,v_n)$. 
First, for $i=1,\dots,n$ define
\[ u'_i = u_i - \sum_{j=1}^{i-1} \frac{\langle u_i,u'_j\rangle_J}{\langle 
u'_j,u'_j\rangle_J} u'_j, \]
where $\langle u_i, u_j \rangle_J = g_J(u_i,u_j)$.
Then, 
\[ \mathcal G_J(f^{\mathrm{sp}}) = \biggl[ 
\frac{u'_1}{|u'_1|_J}, \dots, \frac{u'_n}{|u'_n|_J},
-\frac{J u'_1}{|u'_1|_J}, \dots, -\frac{J u'_n}{|u'_n|_J}
\biggr],
\]
where $|u'_i|_J^2 = \langle u'_i, u'_i \rangle_J$.
\end{rem}

We can now prove the following.

\begin{prop}
The map $\kappa_J : \Lambda_J \to \Lambda_{pl}$, defined by 
\[
\kappa_{J}([u_1, \dots, u_n, v_1, \dots, v_n]) = \mathrm{span}\{u_1, \dots, u_n \},
\]
is a homeomorphism.
\end{prop}

\begin{proof}
Denote by $q_J$ the quotient map from $Fr_J^u$ to $\Lambda_J$. 
Since $\mathcal G_J$ and $q_J$ are surjective submersions, $q_J \circ \mathcal G_J : Fr^{\mathrm{sp}} \to \Lambda_J$ is a quotient map. 
The quotient maps $q_J \circ \mathcal G_J$ and $\Pi$ are constant on each other's fibers. 
Therefore, there is a unique homeomorphism $\kappa_J : \Lambda_J \to \Lambda_{pl}$ such that $\Pi = \kappa_J \circ (q_J \circ \mathcal G_J)$,
and it can be directly checked that the given map $\kappa_J$ satisfies the last relation.
\end{proof}

\section{The structure of \texorpdfstring{$Fr^{u}_J$}{the unitary frame bundle}}
\label{sec:structure-unitary-frame-bundle}

Given a symplectic manifold $(M,\omega)$, one can construct a compatible triple $(\omega, J, g_J)$ as follows \cite{Audin2004, CannasDaSilva2001}. First, choose a Riemannian metric $g$ on $M$. This induces a bundle isomorphism $\mathcal{A}$ on $TM$ by $\omega(\cdot, *) = g(\mathcal{A} \cdot, *)$. Then define the almost complex structure by $J = \mathcal{A}^{-1} \sqrt{-\mathcal{A}^2}$ and the compatible Riemannian metric by $g_J(\dot, *) = \omega(J \cdot, *)$. Note that, in general, $g_J$ differs from $g$.

In principle, starting with different Riemannian structures yields different compatible almost complex structure $J$ and Riemannian metric $g_J$. Additionally, the smoothness of $J$, as well as the smooth structure of $Fr^u_J$, relies on the smoothness of $\sqrt{-\mathcal{A}^2}$.
In this appendix, we show that $J$ is a smooth structure, and that the structure of the unitary frame bundle $Fr^{u}_J$ is independent of the choice of Riemannian metric $g$.

\subsection{Smoothness of \texorpdfstring{$J$}{J}}
\label{subsec:smoothness-of-J}

The smoothness of $\mathcal{A}$ and $\mathcal{A}^{-1}$ follows directly from the smoothness of $\omega$ and $g$. 
Since $-\mathcal{A}^2$ is smooth and symmetric with respect to $g$, it suffices to establish the smoothness of taking the ``square root'' of a symmetric positive definite matrix $\mathbf{A}$, i.e., the map $\mathbf{A} \mapsto \sqrt{\mathbf{A}}$ with $\mathbf{A} = (\sqrt{\mathbf{A}})^2$.

Denote by $\mathrm{PD}(n)$ the set of all $n\times n$ symmetric positive definite matrices. 
The set $\mathrm{PD}(n)$ is open in $\mathbb{R}^{n(n+1)/2}$. For the smoothness of the splitting, it suffices to show that $\mathfrak{sq} : PD(n) \to PD(n)$, given by $\mathfrak{sq}(\mathbf{P}) = \mathbf{P}^2$, is bijective, and thus $\mathbf{A} \mapsto \sqrt{\mathbf{A}}$ is uniquely defined as the inverse of $\mathfrak{sq}$, and also that $\mathfrak{sq}$ is a local diffeomorphism, that is, its tangent map is nondegenerate everywhere.

The bijectivity of $\mathfrak{sq}$ is a basic fact in linear algebra, which follows from the fact that given any $\mathbf{A}, \mathbf{P} \in \mathrm{PD}(n)$, the identity $\mathbf{A} = \mathbf{P}^2$ implies that $\mathbf{A}$ and $\mathbf{P}$ have identical eigenspaces with the corresponding eigenvalues being $\lambda$ and $\lambda^2$, respectively.

For the nondegeneracy of the tangent map of $\mathfrak{sq}$, first note that for any $\mathbf{A}$, there exists $\mathbf{C} \in \mathbb{SO}(n)$, such that
\[
\mathbf{A} = \mathbf{C}^t 
\begin{bmatrix} \lambda_1 & & \\ & \ddots & \\ & & \lambda_n \end{bmatrix}
\mathbf{C}
=\mathbf{C}^t [\lambda] \mathbf{C},
\
\text{and}
\  
\mathbf{A}^2 = \mathbf{C}^t 
\begin{bmatrix} \lambda_1^2 & & \\ & \ddots & \\ & & \lambda_n^2 \end{bmatrix}
\mathbf{C}.
\]
Then the discussion is concluded by checking that, with $i,j=1,...,n$, the tangent vectors
\[
\frac{d}{dt}\bigg|_{t=0}\big([\lambda]+t[ij]\big)^2 = [\lambda][ij]+[ij][\lambda]
\]
are a linearly independent basis, where, $[ij]$ is the symmetric matrix with $a_{ij}=a_{ji}=1$ while all other entries are zero.

\subsection{Isomorphism of unitary frame bundles}
\label{Bundle $Fr^u_J$}

We want to show that, given any two almost complex structures $J_0,J_1$ on a symplectic manifold $(M,\omega)$, the corresponding unitary frame bundles $Fr^u_{J_0}$ and $Fr^u_{J_1}$ are isomorphic as principal $\mathbb{U}(n)$ bundles.
To this end, we construct a $\mathbb{U}(n)$ bundle $\Xi: \mathfrak{U} \to M\times[0,1]$ such that $Fr^u_{J_i} \cong \mathfrak{U}|_{M\times\{i\}}$ for $i=0,1$,
and then the isomorphism of the two bundles follows as direct consequence, see \cite[Theorem~9.6]{Husemoller1994}.

For this purpose, consider the symplectic vector bundle 
\[ \pi_{\mathfrak{sp}} : \mathfrak{V} \cong TM \times [0,1] \to M \times [0,1], \]
with symplectic form $\omega_p$ on each of the fibers $\pi_{\mathfrak{sp}}^{-1}(p,t) \cong T_pM$.
Additionally, consider the corresponding symplectic frame bundle
\[
\Pi_{\mathfrak{sp}} : Fr^{\mathfrak{sp}} \cong Fr^{\mathrm{sp}} \times [0,1] \to M \times [0,1].
\]
The smoothly parameterized Riemannian metrics
\[ g_t = (1-t) g_{J_0} + t g_{J_1},\ t \in [0,1], \]
then define a smooth Euclidean metric on the vector bundle $\mathfrak{V}$.
Using the construction described at the beginning of this appendix, the metrics $g_t$ yield a compatible triple $(\omega, \hat J_t, \hat g_t)$. Due to the compatibility of $g_0 = g_{J_0}$ with $\omega$ and $J_0$, it holds that $\hat g_0 = g_0$ and $\hat J_0 = J_0$. Similarly, $\hat g_1 = g_1$ and $\hat J_1 = J_1$.
Summarizing,
\begin{equation}
\label{eq:barJ_i=J_i}
\hat{J} |_{TM\times\{i\}} = J_i, \ \text{for $i=0,1$}.
\end{equation}
The principal $\mathbb{U}(n)$ bundle $\mathfrak{U}$ is then defined as the unitary frame bundle of $\mathfrak{V}$ with respect to the triple $(\omega, \hat J_t, \hat g_t)$,
and due to Eq.~\eqref{eq:barJ_i=J_i}, the relation $Fr^u_{J_i} \cong \mathfrak{U} |_{M\times\{i\}}$ for $i=0,1$ is satisfied.

\paragraph*{Data Availability Statement.} This work  proceeds within a theoretical and mathematical approach and we do not generate or analyze any datasets.

\bibliographystyle{plain}
\bibliography{main.bib}

\end{document}